# TIME DISCRETIZATION AND MARKOVIAN ITERATION FOR COUPLED FBSDES

BY CHRISTIAN BENDER AND JIANFENG ZHANG[1]

*Technische Universität Braunschweig and
University of Southern California*

In this paper we lay the foundation for a numerical algorithm to simulate high-dimensional coupled FBSDEs under weak coupling or monotonicity conditions. In particular, we prove convergence of a time discretization and a Markovian iteration. The iteration differs from standard Picard iterations for FBSDEs in that the dimension of the underlying Markovian process does not increase with the number of iterations. This feature seems to be indispensable for an efficient iterative scheme from a numerical point of view. We finally suggest a fully explicit numerical algorithm and present some numerical examples with up to 10-dimensional state space.

**1. Introduction.** Motivated by the aim to simulate high-dimensional coupled forward–backward stochastic differential equations (FBSDEs) we study a time discretization and a Markovian iteration for equations of the form

$$(1.1) \quad \begin{cases} X_t = x + \int_0^t b(s, X_s, Y_s)\, ds + \int_0^t \sigma(s, X_s, Y_s)\, dW_s, \\ Y_t = g(X_T) + \int_t^T f(s, X_s, Y_s, Z_s)\, ds - \int_t^T Z_s\, dW_s, \end{cases}$$

where $b, \sigma, f, g$ are deterministic and Lipschitz continuous functions of linear growth which are additionally supposed to satisfy some weak coupling or monotonicity condition. The solution consists of a triplet $(X, Y, Z)$ of adapted processes, which are called the "forward part," the "backward part" and the "control part" respectively. The presence of the control part $Z$ is crucial to find a nonanticipative solution. It often has an intuitive interpretation, for example, as an investment strategy in financial applications; see

Received August 2006; revised April 2007.
[1]Supported in part by NSF Grant DMS-04-03575.

El Karoui, Peng and Quenez [12]. Note that (1.1) is not in its most general form, since $Z$ does not couple into the forward SDE.

Most of the numerical algorithms for coupled FBSDEs, with the notable exception of Delarue and Menozzi [10], exploit the relation to quasi-linear parabolic PDEs via the Ma–Protter–Yong four-step-scheme [17]. Under appropriate conditions, $(X, Y, Z)$ are connected by

(1.2) $\quad Y_t = u(t, X_t), \qquad Z_t = v(t, X_t) \triangleq u_x(t, X_t)\sigma(t, X_t, u(t, X_t))$,

where $u$ is a classical solution of the PDE

(1.3) $\quad \begin{cases} u_t + \frac{1}{2}\mathrm{trace}(\sigma\sigma^*(t,x,u)u_{xx}) \\ \qquad + u_x b(t,x,u) + f(t,x,u,u_x\sigma(t,x,u)) = 0, \\ u(T,x) = g(x). \end{cases}$

The main focus in these approaches is on the numerical solution of the PDE (1.3); see Douglas, Ma and Protter [11], Milstein and Tretyakov [20] and Ma, Shen and Zhao [18]. Since the PDE approach requires existence of a classical solution to (1.3), there is typically need for some smoothness, boundedness and regularity conditions, such as uniform ellipticity of the differential operator. For low-dimensional problems, under such regularity conditions, the PDE approach may generally be regarded as superior to Monte Carlo simulation concerning accuracy and speed. However, solving (1.3) numerically by standard PDE techniques becomes more difficult, if not impossible, with increasing spatial dimension. Hence, it seems necessary to tackle the FBSDE (1.1) directly by probabilistic means in order to solve (1.1) numerically in situations which are beyond the limitations of the PDE approach.

A natural time discretization of equation (1.1) is

(1.4) $\quad \begin{cases} X_0^n \triangleq x, \\ X_{i+1}^n \triangleq X_i^n + b(t_i, X_i^n, Y_i^n)h + \sigma(t_i, X_i^n, Y_i^n)\Delta W_{i+1}, \\ Y_n^n \triangleq g(X_n^n), \\ \hat{Z}_i^n \triangleq \dfrac{1}{h} E_{t_i}\{Y_{i+1}^n \Delta W_{i+1}\}, \\ Y_i^n \triangleq E_{t_i}\{Y_{i+1}^n + f(t_i, X_i^n, Y_{i+1}^n, \hat{Z}_i^n)h\}, \end{cases}$

where $h \triangleq \frac{T}{n}$ and $t_i \triangleq ih, i = 0, 1, \ldots, n$, and $\Delta W_{i+1} \triangleq W_{t_{i+1}} - W_{t_i}$. Here, of course, $E_{t_i}$ denotes the conditional expectation $E\{\cdot|\mathcal{F}_{t_i}\}$. This time discretization was investigated in detail by Zhang [25] for decoupled FBSDEs. Note that $Z$ in (1.2) and $\hat{Z}^n$ in (1.4) may be considered analogous to each other. Indeed, the expression for $Z$ can be rewritten as the Malliavin derivative $D_t Y_t$ of $Y$ and, applying integration by parts under the conditional expectation,

$$\hat{Z}_i^n = \frac{1}{h} E_{t_i}\left\{ \int_{t_i}^{t_{i+1}} D_t Y_{t_{i+1}}^n \, dt \right\}$$

is a natural discretization of the Malliavin derivative of $Y^n$. Concerning the iteration (1.4), it is crucial to notice that $X$ is discretized forwardly and $Y$ is discretized backwardly. Hence, (1.4) is by no means an explicit discretization in the present situation due to the coupling and therefore cannot be implemented directly. Note, however, that one can rewrite

$$(1.5) \qquad Y_i^n = u_i^n(X_i^n), \qquad \hat{Z}_i^n = v_i^n(X_i^n),$$

where

$$(1.6) \qquad \begin{cases} u_n^n(x) \triangleq g(x), \\ X_{i+1}^{n,i,x} \triangleq x + b(t_i, x, u_i^n(x))h + \sigma(t_i, x, u_i^n(x))\Delta W_{i+1}, \\ Y_{i+1}^{n,i,x} \triangleq u_{i+1}^n(X_{i+1}^{n,i,x}), \\ v_i^n(x) \triangleq \frac{1}{h} E\{Y_{i+1}^{n,i,x} \Delta W_{i+1}\}, \\ u_i^n(x) \triangleq E\{Y_{i+1}^{n,i,x} + f(t_i, x, Y_{i+1}^{n,i,x}, v_i^n(x))h\}. \end{cases}$$

Equation (1.6) is still implicit in $u_i^n$, but truly backward in time. Combined with a local updating technique, it serves as starting point for the probabilistic scheme in Delarue and Menozzi [10]. This type of scheme requires, however, apart from estimating the expectations, a discretization of the state space. Such space discretization may again become prohibitive, when the dimension increases.

We, hence, propose to combine the time discretization (1.4) with an iterative scheme. It is known from results by Antonelli [1] and Pardoux and Tang [22] that, under weak coupling or monotonicity conditions, (1.1) has a unique solution $(X, Y, Z)$ which can be constructed via a Picard iteration

$$(1.7) \qquad \begin{cases} \check{X}_t^m = x + \int_0^t b(s, \check{X}_s^m, \check{Y}_s^{m-1})\,ds + \int_0^t \sigma(s, \check{X}_s^m, \check{Y}_s^{m-1})\,dW_s, \\ \check{Y}_t^m = g(\check{X}_T^m) + \int_t^T f(s, \check{X}_s^m, \check{Y}_s^m, \check{Z}_s^m)\,ds - \int_t^T \check{Z}_s^m\,dW_s, \end{cases}$$

starting at $\check{Y}^0 = 0$. The drawback of (1.7) is that the dimension of the underlying Markovian process increases with the number of iterations $m$. Precisely, one can easily see that $\check{Y}^1$ is a function of $\check{X}^1$ and, hence, the right-hand side of the SDE for $\check{X}^2$ depends on $\check{X}^1$ (through $\check{Y}^1$) and $\check{X}^2$. Proceeding this way, one observes that $\check{X}^m$ generally is not Markovian, but only the extended system $(\check{X}^1, \ldots, \check{X}^m)$ is. Consequently, $\check{Y}_t^m$ is a function $\check{u}^m$ of time and $(\check{X}^1, \ldots, \check{X}^m)$, and therefore, the computational effort to estimate $\check{u}^m$ rapidly increases with the number of Picard iterations. This renders a combination of (1.4) with a Picard iteration like (1.7), which was recently suggested by Riviere [23] in theory, impractical from a numerical point of view. The stochastic control approach in Cvitanić and Zhang [8], which iterates over $Z$, faces the same kind of difficulty.

In this paper we introduce an alternative iteration in a way that the dimension of the underlying Markovian process does not change in the number of iterations. It reads, in discretized form, $u_i^{n,0}(x) = 0$, and

$$(1.8) \quad \begin{cases} X_0^{n,m} \triangleq x, \\ X_{i+1}^{n,m} \triangleq X_i^{n,m} + b(t_i, X_i^{n,m}, u_i^{n,m-1}(X_i^{n,m}))h \\ \qquad\qquad + \sigma(t_i, X_i^{n,m}, u_i^{n,m-1}(X_i^{n,m}))\Delta W_{i+1}, \\ Y_n^{n,m} \triangleq g(X_n^{n,m}), \\ \hat{Z}_i^{n,m} \triangleq \dfrac{1}{h} E_{t_i}\{Y_{i+1}^{n,m}\Delta W_{i+1}\}, \\ Y_i^{n,m} \triangleq E_{t_i}\{Y_{i+1}^{n,m} + f(t_i, X_i^{n,m}, Y_{i+1}^{n,m}, \hat{Z}_i^{n,m})h\}, \\ u_i^{n,m}(X_i^{n,m}) = Y_i^{n,m}. \end{cases}$$

The main advantage is that here $Y_i^{n,m}$ is a function of time and $X_i^{n,m}$, but does not depend on $(X_i^{n,\mu}, \mu = 1, \ldots, m-1)$. Establishing the convergence of this new "Markovian" iteration turns out to be more involved than for the standard Picard iteration, because controlling the Lipschitz constant and the linear growth of $u_i^{n,m}(x)$ uniformly in $i, n, m$ becomes crucial. This is indeed the reason why we cannot allow $Z$ to couple in the forward SDE at the current state of our research.

We also indicate how this discretized Markovian iteration may be transformed into a viable numerical scheme, replacing the conditional expectations by simulation based least squares regression and estimating $u^{n,m}$ this way. Such an estimator was introduced by Carrière [6], Longstaff and Schwartz [16] and Clement, Lamberton and Protter [7] in the context of American options and is applied by Gobet, Lemor and Warin [14] and Bender and Denk [3] for decoupled FBSDEs. Although a convergence analysis for this estimator in the present context of a coupled FBSDE is beyond the scope of this paper, we illustrate by some examples with up to 10-dimensional state space that the proposed numerical algorithm works in practice.

The paper is organized as follows: In Section 2 we state the main results on convergence of the discretized Markovian iteration. The proof is given in several steps in Sections 3–5, where we establish the control of the Lipschitz constant, of the linear growth and the convergence of $u^{n,m}$ to $u^n$ respectively. In Section 6 we investigate the error due to the time discretization. To the best of our knowledge, our convergence theorem is the first of this type for coupled FBSDEs which also holds for a degenerate diffusion coefficient $\sigma$. In Section 7 we spell out the proposed numerical scheme and present some numerical examples in Section 8.

**2. Notation and main results.** The main results of this paper estimate the error of the discretized Markovian iteration (1.8) as the number of time steps $n$ and the number of iterations $m$ tend to infinity. Before we can state these results, we need to fix some notation and discuss some assumptions. From now on we suppose, in the theoretical part, that all processes are one-dimensional. This is only to ease the notation and the attentive reader will easily see that all results hold true for the multi-dimensional case as well. The augmented filtration generated by the Brownian motion is denoted by $\mathbf{F} = \{\mathcal{F}_t, 0 \le t \le T\}$.

The first assumption concerns the Lipschitz continuity and monotonicity of the coefficients. It will be in force throughout the whole paper without further notice. Denote

$$\Delta x \triangleq x_1 - x_2, \qquad \Delta y \triangleq y_1 - y_2, \qquad \Delta z \triangleq z_1 - z_2.$$

ASSUMPTION 2.1. (i) There exist (possibly negative) constants $k_b, k_f$ such that

$$[b(t, x_1, y) - b(t, x_2, y)]\Delta x \le k_b |\Delta x|^2,$$
$$[f(t, x, y_1, z) - f(t, x, y_2, z)]\Delta y \le k_f |\Delta y|^2.$$

(ii) $b, \sigma, f, g$ are uniformly Lipschitz continuous with respect to $(x, y, z)$. In particular, there are constants $K$, $b_y$, $\sigma_x$, $\sigma_y$, $f_x$, $f_z$ and $g_x$ such that

$$|b(t, x_1, y_1) - b(t, x_2, y_2)|^2 \le K|\Delta x|^2 + b_y |\Delta y|^2,$$
$$|\sigma(t, x_1, y_1) - \sigma(t, x_2, y_2)|^2 \le \sigma_x |\Delta x|^2 + \sigma_y |\Delta y|^2,$$
$$|f(t, x_1, y_1, z_1) - f(t, x_2, y_2, z_2)|^2 \le f_x |\Delta x|^2 + K|\Delta y|^2 + f_z |\Delta z|^2,$$
$$|g(x_1) - g(x_2)|^2 \le g_x |\Delta x|^2.$$

(iii) $b(t, 0, 0)$, $\sigma(t, 0, 0)$, $f(t, 0, 0, 0)$ are bounded. In particular, there are constants $b_0$, $\sigma_0$, $f_0$ and $g_0$ such that

$$|b(t, x, y)|^2 \le b_0 + K|x|^2 + b_y |y|^2,$$
$$|\sigma(t, x, y)|^2 \le \sigma_0 + \sigma_x |x|^2 + \sigma_y |y|^2,$$
$$|f(t, x, y, z)|^2 \le f_0 + f_x |x|^2 + K|y|^2 + f_z |z|^2,$$
$$|g(x)|^2 \le g_0 + g_x |x|^2.$$

We emphasize that here $b_y$ et al. are constants, not partial derivatives. Indeed, we will not assume any differentiability conditions throughout this paper. For convenience, we also suppose that $K$ is an upper bound for all the constants above.

For results concerning the error due to the time discretization, we require the following assumption.

ASSUMPTION 2.2. *The coefficients $(b, \sigma, f)$ are uniformly Hölder-$\frac{1}{2}$ continuous with respect to $t$.*

If Assumption 2.2 is in force, we use the same constant $K$ to denote an upper bound of the square of the Hölder constants.

To ensure that a solution of (1.1) exists and the iteration converges, we further impose conditions which guarantee that we are in one of the following five cases:

1. *Small time duration*, that is, $T$ is small.
2. *Weak coupling of $Y$ into the forward SDE*, that is, $b_y$ and $\sigma_y$ are small. In particular, if $b_y = \sigma_y = 0$, then the forward equation in (1.1) does not depend on the backward one and, thus, (1.1) is decoupled.
3. *Weak coupling of $X$ into the backward SDE*, that is, $f_x$ and $g_x$ are small. In particular, if $f_x = g_x = 0$, then the backward equation in (1.1) does not depend on the forward one and, thus, (1.1) is also decoupled. In fact, in this case $Z = 0$ and (1.1) reduces to a decoupled system of a forward SDE and an ODE.
4. *$f$ is strongly decreasing in $y$*, that is, $k_f$ is very negative.
5. *$b$ is strongly decreasing in $x$*, that is, $k_b$ is very negative.

The above conditions will be made precise later.

REMARK 2.1. We emphasize that Assumptions 2.1 and 2.2 alone are not sufficient to guarantee existence of a solution to the FBSDE (1.1). For instance, the one-dimensional linear FBSDE,

$$dX_t = Y_t dt, \qquad dY_t = -X_t\, dt + Z_t\, dW_t, \qquad 0 \leq t \leq \frac{3\pi}{4},$$

$$X_0 = x \neq 0, \qquad Y_{3\pi/4} = -X_{3\pi/4},$$

does not admit a solution; see Ma and Yong [19].

We next present two examples from finance, in which we expect one of the conditions 1–5 to hold. For more details on both examples, we refer to Ma and Yong [19], Chapters 8.3 and 8.4.

EXAMPLE 2.1. (i) Consol rate models: In interest rate modeling the term structure can be determined by a so-called short rate $r$; see, for example, Brigo and Mercurio [5]. This short rate may depend on a long rate $Y^{-1}$ (also called consol rate), which in turn is influenced by the short rate via the relation

$$Y_t = E_t\left\{\int_t^T \exp\left(-\int_t^s r_u\, du\right) ds\right\}.$$

This problem can be cast into the FBSDE framework as
$$dr_t = b(t, r_t, Y_t)\, dt + \sigma(t, r_t, Y_t)\, dW_t,$$
$$dY_t = (r_t Y_t - 1)\, dt + Z_t\, dW_t,$$
$$r_0 = R, \qquad Y_T = 0.$$
In generalization of the Hull–White model [15] one can choose
$$b(t, r, y) = \kappa(\theta(t, y) - r), \qquad \sigma(t, r, y) = \eta(t, y) > 0.$$
Here the drift $b$ contains a mean reverting force $\kappa > 0$, which pushes the short rate to a level which may depend on the long rate. The mean reverting force implies that $b$ is monotonically decreasing in $r$, and so we expect condition 5 to hold.

(ii) Stock coupled with an option: In this example the forward SDE describes a system of stock prices, and the backward SDE describes the price of an option on the stocks, leading in a complete market to the FBSDE
$$dS_t = b(t, Y_t) S_t\, dt + \sigma(t, Y_t) S_t\, dW_t,$$
$$dY_t = (rY_t + Z_t \theta(t, Y_t))\, dt + Z_t\, dW_t,$$
$$S_0 = s, \qquad Y_T = g(S_T),$$
where $r$ denotes the riskless interest rate, $\theta$ the premium of risk, and the drift $b$ and volatility $\sigma$ of the stocks are influenced by the price $Y$ of an option on $S$ with pay-off function $g$. As the drift and the volatility of a stock fluctuate only slightly, we can expect that $Y$ weakly couples into the forward part.

Generically, we will derive the following theorems. The first theorem concerns the convergence of the iteration as $m$ tends to infinity.

THEOREM 2.1. *Under Assumption* 2.1, *let one of the conditions* 1–5 *hold true. Then, for sufficiently small* $h$, (1.6) *has an "essentially" unique solution* $u^n$ *with linear growth and there are constants* $C > 0$ *and* $0 < c < 1$ *such that*
$$\max_{0 \leq i \leq n} |u_i^{n,m}(x) - u_i^n(x)|^2 \leq C(|x|^2 + m)c^m,$$
*where* $u^{n,m}$ *is given by* (1.8).

We will see from the proof that the constant $c$, which determines the rate of convergence, depends on the conditions 1–5. Roughly speaking, the stronger the monotonicity (resp. the weaker the coupling, the smaller the time horizon), the smaller one can choose $c$ and, hence, the faster the iteration converges.

Concerning the error due to the time discretization, we obtain the following:

THEOREM 2.2. *Suppose Assumptions* 2.1, 2.2 *hold true, and one of the conditions* 1–5 *is in force. Then equation* (1.3) *admits a viscosity solution* $u(t,x)$ *with linear growth and there is a constant* $C > 0$ *such that, for sufficiently small* $h$,

$$\max_{0 \leq i \leq n} |u_i^n(x) - u(t_i, x)|^2 \leq C(1 + |x|^2)h.$$

Note that the forward part in (1.6) is discretized by an Euler scheme. Hence, the stated convergence of order $1/2$ for the time discretization is the best rate one can hope for.

Combining these two theorems, one can derive the following with a little extra effort:

THEOREM 2.3. *Under the assumptions of Theorem* 2.2 *FBSDE* (1.1) *has a unique solution* $(X, Y, Z)$ *and there are constants* $C > 0$ *and* $0 < c < 1$ *such that, for sufficiently small* $h$,

$$\sup_{1 \leq i \leq n} E \bigg\{ \sup_{t \in [t_{i-1}, t_i]} [|X_t - X_{i-1}^{n,m}|^2 + |Y_t - Y_{i-1}^{n,m}|^2] \bigg\}$$
$$+ \sum_{i=1}^{n} E \bigg\{ \int_{t_{i-1}}^{t_i} |Z_t - \hat{Z}_{i-1}^{n,m}|^2 \, dt \bigg\}$$
$$\leq C(1 + |x|^2)[mc^m + h].$$

These generic results will be made precise in Theorems 5.1, 6.3 and 6.5 below. We emphasize that none of the above theorems requires nondegeneracy of $\sigma$ and, in principle, $X$ and $W$ can have different dimensions. Moreover, we do not suppose any smoothness or boundedness conditions. However, we also underline again that FBSDE (1.1) does not allow coupling through the control part $Z$.

The proof of convergence for the Markovian iteration, which will be given in Sections 3–5, is rather technical. We therefore briefly outline its proof.

**Strategy of proof**. In a standard Picard iteration, like (1.7), one estimates $|\check{Y}^{m+1} - \check{Y}^m|$ in terms of $|\check{X}^{m+1} - \check{X}^m|$ and then $|\check{X}^{m+1} - \check{X}^m|$ in terms of $|\check{Y}^m - \check{Y}^{m-1}|$. However, applying similar techniques to (1.8) yields only estimates of $|X^{n,m+1} - X^{n,m}|$ in terms of $|u^{n,m}(X^{n,m+1}) - u^{n,m-1}(X^{n,m})|$. Since $Y^{n,m} = u^{n,m}(X^{n,m})$, it seems unavoidable to control the Lipschitz constant of $u^{n,m}$ to obtain estimates in terms of $|Y^{n,m} - Y^{n,m-1}|$.

More precisely, suppose, for the moment, that the step size $n$ is fixed, and $u_i^{n,m}(x)$ are bounded functions for all $m$ and $i$. Then one can derive estimates of $|X_i^{n,m+1} - X_i^{n,m}|$ in terms of $\sup_x |u_i^{n,m}(x) - u_i^{n,m-1}(x)|$ and a uniform (in time) Lipschitz constant $L(u^{n,m})$ of $u^{n,m}$ by applying Lemma

3.2 below. Combining this with estimates for $\sup_x |u_i^{n,m+1}(x) - u_i^{n,m}(x)|$ in terms of $|X_i^{n,m+1} - X_i^{n,m}|$ (derived from Lemma 3.3 below), one obtains

$$\sup_x |u_i^{n,m+1}(x) - u_i^{n,m}(x)|^2 \tag{2.1}$$
$$\leq c(L(u^{n,m})) \sup_x |u_i^{n,m}(x) - u_i^{n,m-1}(x)|^2$$

for some constant $c(L(u^{n,m}))$ which depends on the coefficients of the equation and the Lipschitz constant of $u^{n,m}$. Now we wish to iterate the above estimate. To this end, we need a uniform control $L$ of $L(u^{n,m})$ and conditions on the coefficients which ensure that $c(L) < 1$. Section 3 is devoted to deriving such uniform control of the Lipschitz constants.

In general, the conditions imposed in this paper do not guarantee that $u^{n,m}$ is bounded, and therefore, the use of the sup-norm becomes meaningless. However, one can easily see that $u^{n,m}$ is of linear growth. Given linear growing functions $\varphi_i$, there are constants $G(\varphi)$ and $H(\varphi)$ such that

$$|\varphi_i(x)|^2 \leq G(\varphi)|x|^2 + H(\varphi) \qquad \forall (i,x).$$

We consider linearly growing functions $\varphi^1$ and $\varphi^2$ close to each other, if we can choose $G(\varphi^1 - \varphi^2)$ and $H(\varphi^1 - \varphi^2)$ small. Following similar, but slightly more intricate considerations than the ones leading to (2.1), we can estimate $G(u^{n,m+1} - u^{n,m})$ and $H(u^{n,m+1} - u^{n,m})$ in terms of $G(u^{n,m} - u^{n,m-1})$ and $H(u^{n,m} - u^{n,m-1})$; see Theorem 5.2 below. However, the constant, which replaces the above $c(L(u^{n,m}))$ in these estimates, depends on the Lipschitz constant of $u^{n,m}$ and additionally on the linear growth of $u^{n,m-1}$ through $G(u^{n,m-1})$ and $H(u^{n,m-1})$. Hence, for the general case a uniform control for the linear growth of $u^{n,m}$ is required as well. Such control will be given in Section 4. Then, iterating the above estimates yields the convergence of the Markovian iteration under each of the conditions 1–5, as will be demonstrated in Section 5.

In order to study the behavior of the functions $u^{n,m}$, as outlined above, we introduce an important operator $F_n$ for each $n$. For any measurable functions $\varphi = \{\varphi_i\}_{0 \leq i \leq n-1}$, define $\psi$ and $\Phi$ as follows:

$$\begin{cases} \Phi_n(x) \triangleq g(x), \\ X_{i+1}^{\varphi,i,x} \triangleq x + b(t_i, x, \varphi_i(x))h + \sigma(t_i, x, \varphi_i(x))\Delta W_{i+1}, \\ Y_{i+1}^{\varphi,i,x} \triangleq \Phi_{i+1}(X_{i+1}^{\varphi,i,x}), \\ \psi_i(x) \triangleq \frac{1}{h} E\{Y_{i+1}^{\varphi,i,x} \Delta W_{i+1}\}, \\ \Phi_i(x) \triangleq E\{Y_{i+1}^{\varphi,i,x} + f(t_i, x, Y_{i+1}^{\varphi,i,x}, \psi_i(x))h\}. \end{cases} \tag{2.2}$$

We finally set $F_n(\varphi) \triangleq \Phi$. It is then obvious that $u^{n,m} = F_n(u^{n,m-1})$, and $F_n(u^n) = u^n$ if (1.6) has a solution $u^n$. We also point out that $Y^{n,m}$, given by (1.8), can be expressed in the form

$$(2.3) \quad Y_i^{n,m} = Y_{i+1}^{n,m} + f(t_i, X_i^{n,m}, Y_{i+1}^{n,m}, \hat{Z}_i^{n,m})h - \int_{t_i}^{t_{i+1}} Z_t^{n,m} \, dW_t,$$

thanks to the martingale representation theorem. The analogous expression holds for $Y^n$ defined in (1.4).

**3. Lipschitz continuity.** In this section we obtain a Lipschitz constant of $u_i^{n,m}(x)$, uniformly in $(i, n, m)$. To this end, we first investigate the Lipschitz continuity of $F_n(\varphi)$. Given Lipschitz continuous $\varphi$, let $L(\varphi_i)$ denote the square of a Lipschitz constant of $\varphi_i$, and $L(\varphi) \triangleq \sup_i L(\varphi_i)$. Denote

$$(3.1) \quad \begin{aligned} L_0 &\triangleq [b_y + \sigma_y][g_x + f_x T] T e^{[b_y+\sigma_y][g_x+f_x T]T + [2k_b+2k_f+2+\sigma_x+f_z]T}, \\ L_1 &\triangleq [g_x + f_x T][e^{[b_y+\sigma_y][g_x+f_x T]T + [2k_b+2k_f+2+\sigma_x+f_z]T+1} \vee 1]. \end{aligned}$$

Our aim is to derive the following theorem:

THEOREM 3.1. *If*

$$(3.2) \quad L_0 < e^{-1},$$

*then for any $\bar{L} > L_1$ and for $h$ small enough, we have*

$$L(u^{n,m}) \leq \bar{L} \qquad \forall m.$$

Notice that (3.2) holds true in all five cases of Section 2.

We prepare the proof of Theorem 3.1 with two lemmas. For constants $\lambda_j > 0, j = 1, 2, 3$, denote

$$(3.3) \quad \begin{aligned} A_1 &\triangleq 2k_b + \sigma_x + 1 + Kh, \\ A_2 &\triangleq b_y + \sigma_y + Kh, \\ A_3 &\triangleq \lambda_2 + \lambda_3 + (1 + \lambda_2^{-1})Kh, \\ A_4 &\triangleq 2k_f + 1 + \lambda_3^{-1} f_z + (1 + \lambda_2^{-1})Kh, \\ A_5 &\triangleq f_x + (1 + \lambda_2^{-1})Kh. \end{aligned}$$

LEMMA 3.2. *Fix $i$ and for $l = 1, 2$, let*

$$X_{i+1}^l \triangleq X_i^l + b(t_i, X_i^l, \varphi^l(X_i^l))h + \sigma(t_i, X_i^l, \varphi^l(X_i^l))\Delta W_{i+1},$$

*where $X_i^l$ is $\mathcal{F}_{t_i}$-measurable. Assume $\varphi^1$ is uniformly Lipschitz continuous. Then for any $\lambda_1 > 0$,*

$$\begin{aligned} E_{t_i}\{|X_{i+1}^1 - X_{i+1}^2|^2\} &\leq [1 + A_1 h + (1 + \lambda_1) A_2 h L(\varphi^1)] |X_i^1 - X_i^2|^2 \\ &\quad + (1 + \lambda_1^{-1}) A_2 h |\varphi^1(X_i^2) - \varphi^2(X_i^2)|^2. \end{aligned}$$

This lemma can be easily proved by some standard estimates and its proof is therefore omitted.

LEMMA 3.3. *Fix $i$ and for $l = 1, 2$, let*

$$Y_i^l = Y_{i+1}^l + f(t_i, X_i^l, Y_{i+1}^l, \hat{Z}_i^l)h - \int_{t_i}^{t_{i+1}} Z_t^l \, dW_t,$$

*where*

$$\hat{Z}_i^l \triangleq \frac{1}{h} E_{t_i}\{Y_{i+1}^l \Delta W_{i+1}\}.$$

*Then for any $\lambda_2, \lambda_3 > 0$,*

$$|\Delta Y_i|^2 + (1 - A_3)h|\Delta \hat{Z}_i|^2 \leq (1 + A_4 h) E_{t_i}\{|\Delta Y_{i+1}|^2\} + A_5 h |\Delta X_i|^2,$$

*where*

$$\Delta X \triangleq X^1 - X^2, \qquad \Delta Y \triangleq Y^1 - Y^2, \qquad \Delta \hat{Z} \triangleq \hat{Z}^1 - \hat{Z}^2.$$

PROOF. Denote

$$\Delta Z \triangleq Z^1 - Z^2, \qquad \Delta f \triangleq f(t_i, X_i^1, Y_{i+1}^1, \hat{Z}_i^1) - f(t_i, X_i^2, Y_{i+1}^2, \hat{Z}_i^2).$$

Then

(3.4) $$\Delta Y_i + \int_{t_i}^{t_{i+1}} \Delta Z_t \, dW_t = \Delta Y_{i+1} + \Delta f h.$$

Squaring both sides and taking conditional expectation, we have

(3.5)
$$|\Delta Y_i|^2 + E_{t_i}\left\{\int_{t_i}^{t_{i+1}} |\Delta Z_t|^2 dt\right\}$$
$$= E_{t_i}\{|\Delta Y_{i+1}|^2 + 2\Delta Y_{i+1}\Delta f h + |\Delta f|^2 h^2\}.$$

Note that

$$E_{t_i}\left\{\int_{t_i}^{t_{i+1}} |\Delta Z_t|^2 \, dt\right\} \geq \frac{1}{h}\left|E_{t_i}\left\{\int_{t_i}^{t_{i+1}} \Delta Z_t \, dt\right\}\right|^2.$$

By (3.4), we have

$$E_{t_i}\left\{\int_{t_i}^{t_{i+1}} \Delta Z_t \, dt\right\} = E_{t_i}\left\{\int_{t_i}^{t_{i+1}} \Delta Z_t \, dW_t \, \Delta W_{i+1}\right\}$$
$$= E_{t_i}\{[\Delta Y_{i+1} + \Delta f h]\Delta W_{i+1}\}$$
$$= h[\Delta \hat{Z}_i + E_{t_i}\{\Delta f \Delta W_{i+1}\}].$$

Hence,
$$E_{t_i}\left\{\int_{t_i}^{t_{i+1}} |\Delta Z_t|^2\, dt\right\}$$
$$\geq h|\Delta \hat{Z}_i|^2 + 2h\Delta \hat{Z}_i E_{t_i}\{\Delta f \Delta W_{i+1}\}$$
$$\geq (1-\lambda_2)h|\Delta \hat{Z}_i|^2 - \lambda_2^{-1} h|E_{t_i}\{\Delta f \Delta W_{i+1}\}|^2$$
$$\geq (1-\lambda_2)h|\Delta \hat{Z}_i|^2 - \lambda_2^{-1} h^2 E_{t_i}\{|\Delta f|^2\}.$$

Thus, (3.5) implies that

(3.6)
$$|\Delta Y_i|^2 + (1-\lambda_2)h|\Delta \hat{Z}_i|^2$$
$$\leq E_{t_i}\{|\Delta Y_{i+1}|^2 + 2\Delta Y_{i+1}\Delta f h + (1+\lambda_2^{-1})h^2|\Delta f|^2\}.$$

By Assumption 2.1(ii), we have

(3.7) $$|\Delta f|^2 \leq K[|\Delta X_i|^2 + |\Delta Y_{i+1}|^2 + |\Delta \hat{Z}_i|^2].$$

Moreover,
$$\Delta f = \Delta f_1 + \Delta f_2 + \Delta f_3,$$
where
$$\Delta f_1 \triangleq f(t_i, X_i^1, Y_{i+1}^1, \hat{Z}_i^1) - f(t_i, X_i^2, Y_{i+1}^1, \hat{Z}_i^1),$$
$$\Delta f_2 \triangleq f(t_i, X_i^2, Y_{i+1}^1, \hat{Z}_i^1) - f(t_i, X_i^2, Y_{i+1}^2, \hat{Z}_i^1),$$
$$\Delta f_3 \triangleq f(t_i, X_i^2, Y_{i+1}^2, \hat{Z}_i^1) - f(t_i, X_i^2, Y_{i+1}^2, \hat{Z}_i^2).$$

Then by Assumption 2.1(i) and (ii), we get

(3.8)
$$2\Delta Y_{i+1}\Delta f = 2\Delta Y_{i+1}\Delta f_1 + 2\Delta Y_{i+1}\Delta f_2 + 2\Delta Y_{i+1}\Delta f_3$$
$$\leq |\Delta Y_{i+1}|^2 + |\Delta f_1|^2 + 2k_f|\Delta Y_{i+1}|^2$$
$$\quad + \lambda_3^{-1} f_z |\Delta Y_{i+1}|^2 + \lambda_3 f_z^{-1}|\Delta f_3|^2$$
$$\leq |\Delta Y_{i+1}|^2 + f_x|\Delta X_i|^2 + 2k_f|\Delta Y_{i+1}|^2$$
$$\quad + \lambda_3^{-1} f_z |\Delta Y_{i+1}|^2 + \lambda_3 |\Delta \hat{Z}_i|^2.$$

Plugging (3.7) and (3.8) into (3.6), the lemma is proved. □

With these lemmas at hand, we can study the Lipschitz continuity of $F_n(\varphi)$ given Lipschitz continuous $\varphi$.

THEOREM 3.4. *For any Lipschitz continuous $\varphi$, we have*
$$L(F_n(\varphi)) \leq [g_x + A_5 T][\exp([A_1 + A_4 + A_1 A_4 h]T + [A_2 + A_2 A_4 h]T L(\varphi)) \vee 1],$$
*where $\lambda_1 = 0$ and $\lambda_2, \lambda_3 > 0$ are chosen such that*

(3.9) $$A_3 \leq 1.$$

PROOF. Recall (2.2). Fix $i$ and $x_1, x_2$. Denote

$$\Delta x \triangleq x_1 - x_2, \qquad \Delta X \triangleq X^{\varphi,i,x_1} - X^{\varphi,i,x_2}, \qquad \Delta Y \triangleq Y^{\varphi,i,x_1} - Y^{\varphi,i,x_2},$$

$$\Delta \Phi_i \triangleq \Phi_i(x_1) - \Phi_i(x_2), \qquad \Delta \psi_i \triangleq \psi_i(x_1) - \psi_i(x_2).$$

We apply Lemmas 3.2 and 3.3, setting $\lambda_1 = 0$, and obtain

$$E\{|\Delta X_{i+1}|^2\} \leq [1 + A_1 h + A_2 h L(\varphi)]|\Delta x|^2,$$

$$|\Delta \Phi_i|^2 + (1 - A_3) h |\Delta \psi_i|^2 \leq (1 + A_4 h) E\{|\Delta Y_{i+1}|^2\} + A_5 h |\Delta x|^2.$$

By (3.9), we have

$$|\Delta \Phi_i|^2 \leq [1 + A_4 h] L(\Phi_{i+1}) E\{|\Delta X_{i+1}|^2\} + A_5 h |\Delta x|^2$$

$$\leq [1 + A_4 h][1 + A_1 h + A_2 h L(\varphi)] L(\Phi_{i+1}) |\Delta x|^2 + A_5 h |\Delta x|^2.$$

Thus,

(3.10) $$L(\Phi_i) \leq [1 + A_4 h][1 + A_1 h + A_2 h L(\varphi)] L(\Phi_{i+1}) + A_5 h$$
$$\triangleq [1 + \tilde{A} h] L(\Phi_{i+1}) + A_5 h \leq [1 + \tilde{A}^+ h] L(\Phi_{i+1}) + A_5 h,$$

where $\tilde{A}^+ \triangleq \tilde{A} \vee 0$ and

(3.11) $$\tilde{A} \triangleq A_1 + A_4 + A_1 A_4 h + [A_2 + A_2 A_4 h] L(\varphi).$$

Note that $L(\Phi_n) = g_x$. Hence, we can apply the discrete Gronwall inequality to (3.10) and get

$$L(\Phi) \leq e^{\tilde{A}^+ T}[g_x + A_5 T] = [g_x + A_5 T][e^{\tilde{A} T} \vee 1],$$

which, combined with (3.11), yields the assertion. $\square$

We are now in position to give the proof of Theorem 3.1.

PROOF OF THEOREM 3.1. First, by induction, one can easily show that $L_m \triangleq L(u^{n,m}) < \infty$ for each $(n, m)$. Due to Theorem 3.4, we have

$$L_m \leq [g_x + A_5 T][\exp([A_1 + A_4 + A_1 A_4 h]T + [A_2 + A_2 A_4 h] T L_{m-1}) \vee 1],$$

for $\lambda_1 = 0$ and any $\lambda_2, \lambda_3 > 0$ satisfying (3.9).

Introducing

$$\tilde{L}_m \triangleq [A_2 + A_2 A_4 h] T L_m,$$

we get

(3.12) $$\tilde{L}_m \leq [A_2 + A_2 A_4 h][g_x + A_5 T]T[e^{[A_1 + A_4 + A_1 A_4 h]T} e^{\tilde{L}_{m-1}} \vee 1]$$
$$\leq [A_2 + A_2 A_4 h][g_x + A_5 T]T[e^{[A_1 + A_4 + A_1 A_4 h]T} e^{\tilde{L}_{m-1}} + 1].$$

Denote

$$L_0(\lambda, h) \triangleq [A_2 + A_2 A_4 h][g_x + A_5 T]Te^{[A_2+A_2A_4h][g_x+A_5T]T+[A_1+A_4+A_1A_4h]T}.$$

Obviously, $\tilde{L}_0 = 0$. If

(3.13) $$L_0(\lambda, h) \leq e^{-1},$$

then, by induction, one can easily show that

$$\tilde{L}_m \leq [A_2 + A_2 A_4 h][g_x + A_5 T]T + 1 \qquad \forall m.$$

We plug this into the right-hand side of (3.12) to obtain

$$\tilde{L}_m \leq [A_2 + A_2 A_4 h][g_x + A_5 T]T[e^{[A_1+A_4+A_1A_4h]T+[A_2+A_2A_4h][g_x+A_5T]T+1} \vee 1].$$

Thus,

(3.14) $$\begin{aligned}L_m &\leq [g_x + A_5 T][e^{[A_1+A_4+A_1A_4h]T+[A_2+A_2A_4h][g_x+A_5T]T+1} \vee 1] \\ &\triangleq L_1(\lambda, h).\end{aligned}$$

So we want to choose $\lambda_2, \lambda_3$ and $h$ which satisfy (3.9) and minimize $L_0(\lambda, h)$. Recall again that $\lambda_1 = 0$. In dependence of $h$ we set, for small $h$,

(3.15) $$\lambda_2(h) \triangleq \sqrt{h}, \qquad \lambda_3(h) \triangleq 1 - [1+K]\sqrt{h} - Kh.$$

Then $A_3 = 1$ and

$$\lim_{h \downarrow 0} L_0(\lambda(h), h) = L_0, \qquad \lim_{h \downarrow 0} L_1(\lambda(h), h) = L_1.$$

Suppose that (3.2) holds true. Then for any $\bar{L} > L_1$, we obtain $L_0(\lambda(h), h) \leq e^{-1}$ and $L_1(\lambda(h), h) \leq \bar{L}$, provided $h$ is small enough. In view of (3.14), the theorem is proved. $\square$

**4. Linear growth.** This section is devoted to studying the linear growth of the functions $u_i^{n,m}(x)$. Given linearly growing functions $\varphi_i$, assume

$$|\varphi_i(x)|^2 \leq G(\varphi_i)|x|^2 + H(\varphi_i) \qquad \forall x,$$

and let

$$G(\varphi) \triangleq \sup_i G(\varphi_i), \qquad H(\varphi) \triangleq \sup_i H(\varphi_i).$$

To state the main result of this section, we first introduce the functions

(4.1) $$\Gamma_0(x) \triangleq \frac{e^x - 1}{x}, \qquad \Gamma_1(x, y) \triangleq \sup_{0 < \theta < 1} \theta e^{\theta x} \Gamma_0(\theta y) \qquad \forall x, y \in \mathbb{R};$$

and for $G > 0$,

$$c_0(G) \triangleq T[g_x\Gamma_1([2k_f + 1 + f_z]T, [2k_b + 1 + \sigma_x]T + [b_y + \sigma_y]GT)$$
$$+ f_xT\Gamma_0([2k_f + 1 + f_z]T)\Gamma_0([2k_b + 1 + \sigma_x]T + [b_y + \sigma_y]GT)];$$

$$c_1(G) \triangleq [b_y + \sigma_y]c_0(G);$$

$$L_2(G) \triangleq e^{[2k_f+1+f_z]^+T}g_0 + f_0T\Gamma_0([2k_f + 1 + f_z]T) + [b_0 + \sigma_0]c_0(G).$$

THEOREM 4.1. *Assume* (3.2) *holds true and*

(4.2) $$c_1(L_1) < 1.$$

*For any* $\bar{G} > L_1$, $c_1(L_1) < c_1 < 1$, $L_2 > L_2(L_1)$, *and for $h$ small enough, we have*

$$G(u^{n,m}) \leq \bar{G}, \qquad H(u^{n,m}) \leq \frac{L_2}{1 - c_1} \qquad \forall m.$$

Notice that

(4.3)
$$\lim_{x \to -\infty} \Gamma_0(x) = 0,$$
$$\lim_{x \to -\infty} \Gamma_1(x, y) = 0,$$
$$\lim_{y \to -\infty} \Gamma_1(x, y) = 0.$$

Hence, (4.2) is satisfied in cases 1–5 of Section 2.

Again we start with some a-priori estimates whose proofs are fairly straightforward and hence omitted. Denote

(4.4)
$$B_1 \triangleq b_0 + \sigma_0 + Kb_0h,$$
$$B_2 \triangleq f_0 + Kf_0h.$$

LEMMA 4.2. *Assume*

$$X_{i+1} = X_i + b(t_i, X_i, \varphi(X_i))h + \sigma(t_i, X_i, \varphi(X_i))\Delta W_{i+1}.$$

*Then,*

$$E_{t_i}\{|X_{i+1}|^2\} \leq [1 + A_1h + A_2hG(\varphi)]|X_i|^2 + [B_1 + A_2H(\varphi)]h.$$

LEMMA 4.3. *Assume*

$$Y_i = Y_{i+1} + f(t_i, X_i, Y_{i+1}, \hat{Z}_i)h - \int_{t_i}^{t_{i+1}} Z_t \, dW_t,$$

*where*

$$\hat{Z}_i = \frac{1}{h}E_{t_i}\{Y_{i+1}\Delta W_{i+1}\}.$$

*Then, for any* $\lambda_2, \lambda_3 > 0$,
$$|Y_i|^2 + (1 - A_3)h|\hat{Z}_i|^2 \leq [1 + A_4h]E_{t_i}\{|Y_{i+1}|^2\} + A_5h|X_i|^2 + B_2h.$$

To derive bounds for the linear growth of $F_n(\varphi)$, we define discrete time versions of $\Gamma_0$ and $\Gamma_1$ by

(4.5)
$$\Gamma_0^i(x) \triangleq \frac{(1+xh)^i - 1}{x},$$
$$\Gamma_1^n(x, y) \triangleq \sup_{0 \leq i \leq n} (1 + xh)^i \Gamma_0^i(y)$$

and discrete time versions of $c_0(G), c_1(G), L_2(G)$ by

(4.6)
$$c_0(\lambda, h, G) \triangleq g_x \Gamma_1^n(A_4, A_1 + A_2 G) + A_5 \Gamma_0^n(A_4) \Gamma_0^n(A_1 + A_2 G),$$
$$c_1(\lambda, h, G) \triangleq A_2 c_0(\lambda, h, G),$$
$$L_2(\lambda, h, G) \triangleq B_1 c_0(\lambda, h, G) + [e^{A_4 T} \vee 1]g_0 + B_2 \Gamma_0^n(A_4).$$

THEOREM 4.4.  *For any linearly growing* $\varphi$,
(4.7)  $G(F_n(\varphi)) \leq [g_x + A_5 T][e^{[A_1 + A_4 + A_1 A_4 h]T + [A_2 + A_2 A_4 h]T G(\varphi)} \vee 1],$
(4.8)  $H(F_n(\varphi)) \leq c_1(\lambda, h, G(\varphi))H(\varphi) + L_2(\lambda, h, G(\varphi)),$
*where* $\lambda_2, \lambda_3 > 0$ *are supposed to fulfill* (3.9).

PROOF. Denote $\Phi \triangleq F_n(\varphi)$. Fix $(i_0, x)$ and define, for $i = i_0, \ldots, n-1$,

(4.9)
$$\begin{cases} X_{i_0} \triangleq x, \\ X_{i+1} \triangleq X_i + b(t_i, X_i, \varphi_i(X_i))h + \sigma(t_i, X_i, \varphi_i(X_i))\Delta W_{i+1}, \\ Y_n \triangleq g(X_n), \\ \hat{Z}_i \triangleq \frac{1}{h} E_{t_i}\{Y_{i+1} \Delta W_{i+1}\}, \\ Y_i \triangleq Y_{i+1} + f(t_i, X_i, Y_{i+1}, \hat{Z}_i)h - \int_{t_i}^{t_{i+1}} Z_t \, dW_t. \end{cases}$$

Obviously $Y_{i_0} = \Phi_{i_0}(x)$. We obtain from Lemma 4.2 that
$$E\{|X_{i+1}|^2\} \leq [1 + A_1 h + A_2 h G(\varphi)]E\{|X_i|^2\} + [B_1 + A_2 H(\varphi)]h.$$
Then

(4.10)
$$E\{|X_i|^2\} \leq [1 + A_1 h + A_2 h G(\varphi)]^{i-i_0} E\{|X_{i_0}|^2\}$$
$$+ [B_1 + A_2 H(\varphi)]h \sum_{j=i_0}^{i-1} [1 + A_1 h + A_2 h G(\varphi)]^{j-i_0}$$
$$= [1 + A_1 h + A_2 h G(\varphi)]^{i-i_0} |x|^2$$
$$+ [B_1 + A_2 H(\varphi)]\Gamma_0^{i-i_0}(A_1 + A_2 G(\varphi)).$$

Next, applying Lemma 4.3 and by (3.9), we have
$$E\{|Y_i|^2\} \leq [1 + A_4 h]E\{|Y_{i+1}|^2\} + A_5 h E\{|X_i|^2\} + B_2 h.$$
Note that
$$|Y_n|^2 \leq g_0 + g_x |X_n|^2.$$
Then
$$|\Phi_{i_0}(x)|^2 = |Y_{i_0}|^2 \leq (1 + A_4 h)^{n-i_0}[g_0 + g_x E\{|X_n|^2\}]$$
$$+ A_5 h \sum_{i=i_0}^{n-1}(1 + A_4 h)^{i-i_0} E\{|X_i|^2\} + B_2 \Gamma_0^{n-i_0}(A_4).$$
This, together with (4.10), implies
$$G(\Phi_{i_0}) \leq (1 + A_4 h)^{n-i_0} g_x [1 + A_1 h + A_2 h G(\varphi)]^{n-i_0}$$
$$+ A_5 h \sum_{i=i_0}^{n-1}(1 + A_4 h)^{i-i_0}[1 + A_1 h + A_2 h G(\varphi)]^{i-i_0},$$
$$H(\Phi_{i_0}) \leq (1 + A_4 h)^{n-i_0} g_0 + B_2 \Gamma_0^{n-i_0}(A_4)$$
$$+ [B_1 + A_2 H(\varphi)]\bigg[g_x(1 + A_4 h)^{n-i_0}\Gamma_0^{n-i_0}(A_1 + A_2 G(\varphi))$$
$$+ A_5 h \sum_{i=i_0}^{n-1}(1 + A_4 h)^{i-i_0}\Gamma_0^{i-i_0}(A_1 + A_2 G(\varphi))\bigg].$$
Note that, for $0 \leq i \leq n$,
$$(1 + xh)^i \leq e^{xT} \vee 1, \qquad \Gamma_0^i(x) \leq \Gamma_0^n(x),$$
(4.11)
$$(1 + xh)^i \Gamma_0^i(y) \leq \Gamma_1^n(x, y).$$
Then
$$G(\Phi_{i_0}) \leq [g_x + A_5 T][e^{[A_1 + A_4 + A_1 A_4 h]T + [A_2 + A_2 A_4 h]T G(\varphi)} \vee 1],$$
$$H(\Phi_{i_0}) \leq [e^{A_4 T} \vee 1]g_0 + B_2 \Gamma_0^n(A_4) + c_0(\lambda, h, G)[B_1 + A_2 H(\varphi)].$$
Since the right-hand side does not depend on $i_0$, the assertion is proved. $\square$

After these preparations we give the proof of Theorem 4.1.

PROOF OF THEOREM 4.1. Denote $G_m \triangleq G(u^{n,m}), H_m \triangleq H(u^{n,m})$. Obviously, $G_0 = H_0 = 0$. We may now conclude from Theorem 4.4 that, under (3.9),

(4.12) $\quad G_m \leq [g_x + A_5 T][e^{[A_1 + A_4 + A_1 A_4 h]T + [A_2 + A_2 A_4 h]T G_{m-1}} \vee 1],$

(4.13) $\quad H_m \leq c_1(\lambda, h, G_{m-1})H_{m-1} + L_2(\lambda, h, G_{m-1}).$

We now choose $\lambda_3(h)$ and $\lambda_4(h)$ as in (3.15) for small $h$. Since (3.2) holds true, for any $\bar{G} > L_1$, we can follow the same arguments as in Theorem 3.1 and get $G(u^{n,m}) \leq \bar{G}$ from (4.12). Note that

$$\lim_{n \to \infty} \Gamma_0^n(x) = T\Gamma_0(xT), \qquad \lim_{n \to \infty} \Gamma_1^n(x,y) = T\Gamma_1(xT, yT),$$
$$\lim_{h \downarrow 0} c_1(\lambda(h), h, G) = c_1(G), \qquad \lim_{h \downarrow 0} L_2(\lambda(h), h, G) = L_2(G).$$

For any $c_1$, $c_1(L_1) < c_1 < 1$, and $L_2$, $L_2(L_1) < L_2$, we can choose $\bar{G} > L_1$ such that $c_1(\bar{G}) < c_1$ and $L_2(\bar{G}) < L_2$. Then, for sufficiently small $h$, it holds that $c_1(\lambda(h), h, \bar{G}) \leq c_1$ and $L_2(\lambda(h), h, \bar{G}) \leq L_2$. Now, by (4.13), we get

$$H_m \leq c_1 H_{m-1} + L_2,$$

which implies the result. $\square$

**5. Convergence of the Markovian iteration.** We now make the assumptions of Theorem 2.1 precise and prove convergence of the Markovian iteration as the number of iteration steps tends to infinity.

To this end, we first introduce

$$c_2(\lambda_1, L, G) \triangleq [e^{[2k_b + 1 + \sigma_x + [b_y + \sigma_y]G]T} \vee 1](1 + \lambda_1^{-1})[b_y + \sigma_y]T$$
$$\times [g_x \Gamma_1([2k_f + 1 + f_z]T,$$
$$[2k_b + 1 + \sigma_x + (1 + \lambda_1)[b_y + \sigma_y]L]T)$$
$$+ f_x T \Gamma_0([2k_f + 1 + f_z]T)$$
$$\times \Gamma_0([2k_b + 1 + \sigma_x + (1 + \lambda_1)[b_y + \sigma_y]L]T)];$$
$$c_2(L, G) \triangleq \inf_{\lambda_1 > 0} c_2(\lambda_1, L, G).$$

We will prove the following theorem:

THEOREM 5.1. *Assume* (3.2) *and*

(5.1) $$c_2(L_1, L_1) < 1.$$

(i) *For any* $\bar{L} > L_1, \bar{G} > L_1, L_2 > L_2(L_1), c_1(L_1) < c_1 < 1$, *there exists a solution* $u^n$ *to* (1.6) *such that*

(5.2) $$L(u^n) \leq \bar{L}, \qquad G(u^n) \leq \bar{G}, \qquad H(u^n) \leq \bar{H} \triangleq \frac{L_2}{1 - c_1},$$

*if $h$ is small enough.*

(ii) *For any $c_2(L_1, L_1) < c_2 < 1$ and for $h$ small enough,*

$$\max_{0 \le i \le n} |u_i^{n,m}(x) - u_i^n(x)|^2$$

(5.3)
$$\le \frac{3\bar{G}}{(1 - \sqrt{c_2})^2} |x|^2 c_2^m$$

$$+ \frac{3}{(1 - \sqrt{c_2})^4} \left[ \frac{\bar{H}}{m} + [b_0 + \sigma_0 + (b_y + \sigma_y)\bar{H}]T\bar{G} \right] mc_2^m.$$

(iii) *Fix $G > 0$ and suppose $\tilde{u}^n$ is another solution to (1.6) with linear growth such that $G(\tilde{u}^n) \le G$. Then $\tilde{u}^n = u^n$, if $h$ (depending on $G$) is small enough.*

REMARK 5.1. (i) In view of (4.3), it is straightforward to see that (5.1) is also satisfied in cases 1–5 of Section 2.
(ii) One can check directly that $c_1(L) \le c_2(L, L)$, and thus, condition (5.1) implies (4.2).

Again we first study the operator $F_n$ to prepare the proof of Theorem 5.1.

THEOREM 5.2. *Assume $\varphi^1, \varphi^2$ have linear growth and $\varphi^1$ is Lipschitz continuous. Then, for any $\lambda_1 > 0$,*

$$G(F_n(\varphi^1) - F_n(\varphi^2))$$
$$\le c_2(\lambda_1, h, L(\varphi^1), G(\varphi^2))G(\varphi^1 - \varphi^2),$$
$$H(F_n(\varphi^1) - F_n(\varphi^2))$$
$$\le c_2(\lambda_1, h, L(\varphi^1), G(\varphi^2))H(\varphi^1 - \varphi^2)$$
$$+ c_2(\lambda_1, h, L(\varphi^1), G(\varphi^2))[B_1 + A_2H(\varphi^2)]TG(\varphi^1 - \varphi^2),$$

*where $\lambda_2, \lambda_3$ are chosen such that (3.9) holds, and*

(5.4)
$$c_2(\lambda_1, h, L, G) \triangleq [e^{[A_1 + A_2G]T} \vee 1](1 + \lambda_1^{-1})A_2$$
$$\times [g_x \Gamma_1^n(A_4, A_1 + (1 + \lambda_1)A_2L)$$
$$+ A_5 \Gamma_0^n(A_4)\Gamma_0^n(A_1 + (1 + \lambda_1)A_2L)].$$

PROOF. For $l = 1, 2$, denote $\Phi^l \triangleq F_n(\varphi^l)$. Fix $(i_0, x)$ and define $(X^l, Y^l, \hat{Z}^l)$ analogously to (4.9). Then obviously $Y_{i_0}^l = \Phi_{i_0}^l(x)$. Denote $L \triangleq L(\varphi_1)$, and

$$\Delta X \triangleq X^1 - X^2, \qquad \Delta Y \triangleq Y^1 - Y^2, \qquad \Delta \hat{Z} \triangleq \hat{Z}^1 - \hat{Z}^1,$$
$$\Delta \varphi \triangleq \varphi^1 - \varphi^2, \qquad \Delta \Phi \triangleq \Phi^1 - \Phi^2.$$

Application of Lemma 3.2 yields

$$
\begin{aligned}
(5.5) \quad & E\{|\Delta X_{i+1}|^2\} \\
& \leq E\{[1 + A_1 h + (1 + \lambda_1)A_2 h L]|\Delta X_i|^2 + (1 + \lambda_1^{-1})A_2 h |\Delta \varphi(X_i^2)|^2\}.
\end{aligned}
$$

Note that

$$|\Delta \varphi(X_i^2)|^2 \leq G(\Delta\varphi)|X_i^2|^2 + H(\Delta\varphi).$$

By the first inequality in (4.10) and (4.11),

$$(5.6) \quad \sup_{i_0 \leq i \leq n} E\{|X_i^2|^2\} \leq [|x|^2 + [B_1 + A_2 H(\varphi^2)]T][e^{[A_1 + A_2 G(\varphi^2)]T} \vee 1] \triangleq \tilde{A}.$$

Then (5.5) implies

$$
\begin{aligned}
E\{|\Delta X_{i+1}|^2\} & \leq [1 + A_1 h + (1 + \lambda_1)A_2 h L] E\{|\Delta X_i|^2\} \\
& \quad + (1 + \lambda_1^{-1})A_2 h [G(\Delta\varphi)\tilde{A} + H(\Delta\varphi)].
\end{aligned}
$$

Since $\Delta X_{i_0} = 0$, we get

$$
\begin{aligned}
(5.7) \quad & \sup_{i_0 \leq i \leq n} E\{|\Delta X_i|^2\} \\
& \leq (1 + \lambda_1^{-1})A_2 h [G(\Delta\varphi)\tilde{A} + H(\Delta\varphi)] \\
& \quad \times \sum_{i=i_0}^{n-1} [1 + A_1 h + (1 + \lambda_1)A_2 h L]^{i-i_0} \\
& = (1 + \lambda_1^{-1})A_2 [G(\Delta\varphi)\tilde{A} + H(\Delta\varphi)]\Gamma_0^{n-i_0}(A_1 + (1 + \lambda_1)A_2 L).
\end{aligned}
$$

Furthermore, we obtain from Lemma 3.3 and (3.9),

$$E\{|\Delta Y_i|^2\} \leq [1 + A_4 h] E\{|\Delta Y_{i+1}|^2\} + A_5 h E\{|\Delta X_i|^2\}.$$

Hence, by (5.7) and (4.11),

$$
\begin{aligned}
|\Delta \Phi_{i_0}(x)|^2 & = |\Delta Y_{i_0}|^2 \\
& \leq (1 + A_4 h)^{n-i_0} E\{|\Delta Y_n|^2\} + A_5 \Gamma_0^{n-i_0}(A_4) \sup_{i_0 \leq i \leq n} E\{|\Delta X_i|^2\} \\
& \leq [(1 + A_4 h)^{n-i_0} g_x + A_5 \Gamma_0^{n-i_0}(A_4)] \sup_{i_0 \leq i \leq n} E\{|\Delta X_i|^2\} \\
& \leq (1 + \lambda_1^{-1})A_2 [g_x(1 + A_4 h)^{n-i_0}\Gamma_0^{n-i_0}(A_1 + (1 + \lambda_1)A_2 L) \\
& \quad + A_5 \Gamma_0^{n-i_0}(A_4)\Gamma_0^{n-i_0}(A_1 + (1 + \lambda_1)A_2 L)] \\
& \quad \times [G(\Delta\varphi)\tilde{A} + H(\Delta\varphi)] \\
& \leq (1 + \lambda_1^{-1})A_2 [G(\Delta\varphi)\tilde{A} + H(\Delta\varphi)] \\
& \quad \times [g_x \Gamma_1^n(A_4, A_1 + (1 + \lambda_1)A_2 L) \\
& \quad + A_5 \Gamma_0^n(A_4)\Gamma_0^n(A_1 + (1 + \lambda_1)A_2 L)],
\end{aligned}
$$

which, together with (5.6), implies the theorem. □

We can apply this theorem to estimate the distance between $u^{n,m}$ and $u^{n,m-1}$.

THEOREM 5.3. *Assume that $L(u^{n,m}) \leq \bar{L}$, $G(u^{n,m}) \leq \bar{G}$ and $H(u^{n,m}) \leq \bar{H}$ for all $m \in \mathbb{N}$ and sufficiently small $h$. Moreover, assume*

$$c_2(\bar{L}, \bar{G}) < 1.$$

*Then for any $c_2(\bar{L}, \bar{G}) < c_2 < 1$ and for $h$ small enough,*

(5.8) $\quad G(u^{n,m+1} - u^{n,m}) \leq \bar{G} c_2^m;$

(5.9) $\quad H(u^{n,m+1} - u^{n,m}) \leq [\bar{H} + [(b_0 + \sigma_0) + (b_y + \sigma_y)\bar{H}]T\bar{G}m]c_2^m.$

PROOF. Choose $\lambda_2, \lambda_3$ depending on $h$ as in (3.15). Note that with this choice

$$\lim_{h \to 0} \Gamma_1^n(x, y) = T\Gamma_1(xT, yT), \qquad \lim_{h \to 0} c_2(\lambda_1, h, L, G) = c_2(\lambda_1, L, G).$$

Hence, we may find an appropriate $\lambda_1$ such that, for $h$ small enough,

$$c_2(\lambda_1, h, \bar{L}, \bar{G}) \leq c_2,$$
$$c_2(\lambda_1, h, \bar{L}, \bar{G})[B_1 + A_2\bar{H}] \leq c_2[[b_0 + \sigma_0] + [b_y + \sigma_y]\bar{H}].$$

Denote

$$\Delta u^{n,m} \triangleq u^{n,m} - u^{n,m-1}.$$

Applying Theorem 5.2, we get, for small $h$,

(5.10) $\quad G(\Delta u^{n,m+1}) \leq c_2 G(\Delta u^{n,m}),$

(5.11) $\quad H(\Delta u^{n,m+1}) \leq c_2 H(\Delta u^{n,m}) + c_2[b_0 + \sigma_0 + [b_y + \sigma_y]\bar{H}]TG(\Delta u^{n,m}).$

Note that

$$G(\Delta u^{n,1}) = G(u^{n,1}) \leq \bar{G},$$
$$H(\Delta u^{n,1}) = H(u^{n,1}) \leq \bar{H}.$$

By (5.10), we get (5.8). Moreover, together with (5.8), (5.11) implies (5.9) immediately. The proof is complete now. □

Theorem 5.1 can now be proved by iterating the above theorem.

PROOF OF THEOREM 5.1. Assume $\bar{G}, \bar{L}, L_2, c_1, c_2$ satisfy the conditions specified in the theorem. Without loss of generality, we assume $c_2(\bar{L}, \bar{G}) < c_2$.

Recall that (5.1) implies (4.2). Hence, by Theorems 3.1 and 4.1, we get, for $h$ small enough,

$$L(u^{n,m}) \leq \bar{L}, \qquad G(u^{n,m}) \leq \bar{G}, \qquad H(u^{n,m}) \leq \bar{H}.$$

Hence, (i) will follow directly from (ii).

To prove (ii), we denote

$$\tilde{L} \triangleq [[b_0 + \sigma_0] + [b_y + \sigma_y]\bar{H}]T\bar{G}.$$

Applying Theorem 5.3, we get

$$|u_i^{n,m+1}(x) - u_i^{n,m}(x)|^2 \leq [\bar{G}|x|^2 + \bar{H} + \tilde{L}m]c_2^m.$$

Then

$$|u_i^{n,m+1}(x) - u_i^{n,m}(x)| \leq [\sqrt{\bar{G}}|x| + \sqrt{\bar{H}} + \sqrt{\tilde{L}m}]c_2^{m/2}.$$

Thus, for any $m_1 > m$,

$$|u_i^{n,m}(x) - u_i^{n,m_1}(x)| \leq \sum_{j=m}^{\infty} \left[\sqrt{\bar{G}}|x| + \sqrt{\bar{H}} + \sqrt{\frac{\tilde{L}}{m}j}\right]c_2^{j/2}$$

$$\leq [\sqrt{\bar{G}}|x| + \bar{H}]\frac{c_2^{m/2}}{1 - \sqrt{c_2}} + \sqrt{\frac{\tilde{L}}{m}}\frac{m(1 - \sqrt{c_2}) + \sqrt{c_2}}{(1 - \sqrt{c_2})^2}c_2^{m/2}.$$

Note that the right-hand side above converges to 0 as $m \to \infty$. Then $u_i^{n,m}(x)$ is a Cauchy sequence and hence converges to some $u_i^n(x)$. Moreover,

$$|u_i^{n,m}(x) - u_i^n(x)|^2 \leq 3[\bar{G}|x|^2 + \bar{H} + \tilde{L}m]\frac{c_2^m}{(1 - \sqrt{c_2})^4},$$

which leads to (5.3) and thus proves (ii).

It remains to prove (iii). For any $G > 0$, assume $\tilde{u}^n$ is another solution to (1.6) with linear growth such that $G(\tilde{u}^n) \leq G$. Then $F_n(\tilde{u}^n) = \tilde{u}^n$. Note that $\tilde{u}_n^n = g = u_n^n$. Assume $\tilde{u}_{i+1}^n = u_{i+1}^n$. We now apply a local version of Theorem 5.2. That is, we consider (2.2) only on the interval $[t_i, t_{i+1}]$ with terminal condition $\Phi_{i+1}(x) \triangleq u_{i+1}^n(x)$ (instead of on $[0, T]$ with terminal condition $g(x)$). We note that in this case there is only one time subinterval. One can check directly that

$$\Gamma_0^1(x) = h, \qquad \Gamma_1^1(x, y) = (1 + xh)h.$$

Then (5.4) becomes

$$\tilde{c}_2(h) \triangleq [e^{[A_1 + A_2 G]h} \vee 1](1 + \lambda_1^{-1})A_2[\bar{L}(1 + A_4 h)h + A_5 h^2].$$

Set $\varphi^1 \triangleq u^n, \varphi^2 \triangleq \tilde{u}^n$. We note that, for given $\varphi$, $G(\varphi)$ and $H(\varphi)$ are not unique. So, in general, Theorem 5.2 does not lead to

$$G(u_i^n - \tilde{u}_i^n) \leq \tilde{c}_2(h) G(u_i^n - \tilde{u}_i^n),$$
$$H(u_i^n - \tilde{u}_i^n) \leq \tilde{c}_2(h) H(u_i^n - \tilde{u}_i^n) + \tilde{c}_2(h)[B_1 + A_2 H(\tilde{u}^n)] T G(u_i^n - \tilde{u}_i^n).$$

To get around this difficulty, we use the definitions of $G(\varphi), H(\varphi)$. Let $G_0, H_0$ be some constants satisfying

$$|u_i^n(x) - \tilde{u}_i^n(x)|^2 \leq G_0 |x|^2 + H_0.$$

For $\nu = 1, 2, \ldots$, denote

$$G_\nu \triangleq \tilde{c}_2(h) G_{\nu-1}, \qquad H_\nu \triangleq \tilde{c}_2(h) H_{\nu-1} + \tilde{c}_2(h)[B_1 + A_2 H(\tilde{u}^n)] T G_{\nu-1}.$$

Now applying Theorem 5.2 repeatedly for $\nu = 1, 2, \ldots$, we get

(5.12) $$|u_i^n(x) - \tilde{u}_i^n(x)|^2 \leq G_\nu |x|^2 + H_\nu \qquad \forall \nu.$$

Note that

$$G_\nu = G_0 \tilde{c}_2(h)^\nu, \qquad H_\nu = H_0 \tilde{c}_2(h)^\nu + [B_1 + A_2 H(\tilde{u}^n)] T \nu \tilde{c}_2(h)^\nu.$$

For $h$ small enough, we have $\tilde{c}_2(h) < 1$. Then

$$\lim_{\nu \to \infty} G_\nu = \lim_{\nu \to \infty} H_\nu = 0.$$

Sending $\nu \to \infty$ in (5.12), we get $\tilde{u}_i^n = u_i^n$. Repeating the arguments backwardly, we prove that $\tilde{u}^n = u^n$. □

**6. Convergence of the time discretization.** We now study the error due to the time discretization. We first introduce a continuous time version of the operator $F_n$. Suppose $\varphi$ is a function on $[0, T] \times \mathbb{R}$ which is Lipschitz in the space variable and let $(X^{\varphi,r,x}, Y^{\varphi,r,x}, Z^{\varphi,r,x})$ be the unique solution to the decoupled FBSDE ($0 \leq r \leq t \leq T$)

$$\begin{cases} X_t^{\varphi,r,x} = x + \int_r^t b(s, X_s^{\varphi,r,x}, \varphi(s, X_s^{\varphi,r,x})) \, ds \\ \qquad + \int_r^t \sigma(s, X_s^{\varphi,r,x}, \varphi(s, X_s^{\varphi,r,x})) \, dW_s, \\ Y_t^{\varphi,r,x} = g(X_T^{\varphi,r,x}) + \int_t^T f(s, X_s^{\varphi,r,x}, Y_s^{\varphi,r,x}, Z_s^{\varphi,r,x}) \, ds - \int_t^T Z_s^{\varphi,r,x} \, dW_s. \end{cases}$$

We then define $\Phi(t, x) \triangleq Y_t^{\varphi,t,x}$ and $F(\varphi) \triangleq \Phi$. It is known from Pardoux and Peng [21] that, under Assumption 2.2 and if $\varphi$ is additionally continuous as a function in time and space, $\Phi$ is a viscosity solution to the following semilinear PDE:

$$\begin{cases} \Phi_t + \frac{1}{2} \sigma^2(t, x, \varphi) \Phi_{xx} + b(t, x, \varphi) \Phi_x + f(t, x, \Phi, \Phi_x \sigma(t, x, \varphi)) = 0, \\ \Phi(T, x) = g(x). \end{cases}$$

We now define recursively $\bar{u}^0 \triangleq 0$ and $\bar{u}^m \triangleq F(\bar{u}^{m-1})$. Then the following theorem can be proved similarly to, actually more easily than, Theorem 5.1. A detailed proof can be found in the appendix of the preprint version, which is available from the authors upon request.

THEOREM 6.1. *Assume* (3.2) *and* (5.1) *hold true.*

(i) $\bar{u}^m$ *converges to some function $u$ uniformly on compacts.*
(ii) $|u(t, x_1) - u(t, x_2)|^2 \leq L_1 |x_1 - x_2|^2$; $|u(t, x)|^2 \leq L_1 |x|^2 + \frac{L_2(L_1)}{1 - c_1(L_1)}$.
(iii) $|u(t, x) - u(s, x)|^2 \leq C(1 + |x|^2)|t - s|$ *for some constant $C$.*
(iv) $F(u) = u$. *Moreover, if $F(\tilde{u}) = \tilde{u}$ and $\tilde{u}$ has linear growth, then $\tilde{u} = u$.*
(v) *Under Assumption* 2.2, $u$ *is a viscosity solution to* (1.3).

From now on we denote by $C$ a generic constant which may depend on the coefficients $b, \sigma, f, g$, but is independent of $n, h$ and $x$. The value of $C$ may vary from line to line.

We next consider the following decoupled FBSDE:

$$
(6.1) \quad \begin{cases} X_t = x + \int_0^t b(s, X_s, u(s, X_s))\, ds + \int_0^t \sigma(s, X_s, u(s, X_s))\, dW_s, \\ Y_t = g(X_T) + \int_t^T f(s, X_s, Y_s, Z_s)\, ds - \int_t^T Z_s\, dW_s, \end{cases}
$$

and its time discretization:

$$
(6.2) \quad \begin{cases} \tilde{X}_0^n \triangleq x, \\ \tilde{X}_{i+1}^n \triangleq \tilde{X}_i^n + b(t_i, \tilde{X}_i^n, u(t_i, \tilde{X}_i^n))h + \sigma(t_i, \tilde{X}_i^n, u(t_i, \tilde{X}_i^n))\Delta W_{i+1}, \\ \tilde{Y}_n^n \triangleq g(\tilde{X}_n^n), \\ \tilde{Z}_i^n \triangleq \frac{1}{h} E_{t_i}\{\tilde{Y}_{i+1}^n \Delta W_{i+1}\}, \\ \tilde{Y}_i^n \triangleq E_{t_i}\{\tilde{Y}_{i+1}^n + f(t_i, \tilde{X}_i^n, \tilde{Y}_{i+1}^n, \tilde{Z}_i^n)h\}. \end{cases}
$$

Denote $u_i^0(x) \triangleq u(t_i, x)$ and $\tilde{u}^n \triangleq F_n(u^0)$. It is obvious that $\tilde{Y}_i^n = \tilde{u}_i^n(\tilde{X}_i^n)$. Note again that (6.1) is decoupled. By Theorem 6.1, and applying nowadays standard arguments for decoupled FBSDEs (see, e.g., Delarue [9] and Zhang [26] for (i), and Zhang [25] and Bouchard and Touzi [4] for (ii)), we can derive the following corollary. A detailed proof is again given in the appendix of the preprint version.

COROLLARY 6.2. *Assume all the conditions in Theorem* 6.1 *hold true.*

(i) *FBSDE* (1.1) *has a unique solution* $(X, Y, Z)$, *which also solves* (6.1), *and it holds that* $Y_t = u(t, X_t)$.
(ii) *Moreover, we have the following estimates:*

$$
(6.3) \quad |\tilde{u}_i^n(x) - u(t_i, x)|^2 \leq C(1 + |x|^2)h,
$$

$$\sup_{1\leq i\leq n} E\bigg\{\sup_{t\in[t_{i-1},t_i]}[|X_t - \tilde{X}^n_{i-1}|^2 + |Y_t - \tilde{Y}^n_{i-1}|^2]\bigg\}$$

(6.4)
$$+ \sum_{i=1}^n E\bigg\{\int_{t_{i-1}}^{t_i} |Z_t - \tilde{Z}^n_{i-1}|^2 dt\bigg\}$$

$$\leq C(1+|x|^2)h.$$

Applying the above decoupling relation, $Y_t = u(t, X_t)$, and the convergence results for decoupled FBSDEs, stated in (6.3)–(6.4), we can establish the convergence of $u^n$, as the time grid becomes finer.

THEOREM 6.3. *Suppose Assumption* 2.2 *is in force, and* (3.2) *and* (5.1) *hold true. Then*

$$|u_i^n(x) - u(t_i, x)|^2 \leq C[1+|x|^2]h.$$

PROOF. For any $\bar{L} > c_1(L_1)$ and $\bar{G} > c_1(L_1)$, when $h$ is small, we have

$$L(u^n) \leq \bar{L}, \qquad G(u^n) \leq \bar{G}.$$

Moreover, we know from Theorem 6.1(ii) that

$$L(u^0) \leq L_1, \qquad G(u^0) \leq L_1, \qquad H(u^0) \leq \bar{H} \triangleq \frac{L_2(L_1)}{1 - c_1(L_1)}.$$

Note that $F_n(u^n) = u^n$. Applying Theorem 5.2 on $u^n$ and $u^0$, we get

$$G(u^n - \tilde{u}^n) \leq c_2(\lambda_1, h, \bar{L}, \bar{G})G(u^n - u^0),$$
$$H(u^n - \tilde{u}^n) \leq c_2(\lambda_1, h, \bar{L}, \bar{G})H(u^n - u^0)$$
$$\qquad + c_2(\lambda_1, h, \bar{L}, \bar{G})[B_1 + A_2\bar{H}]TG(u^n - u^0).$$

For any $\varepsilon > 0$, we obtain, thanks to (6.3),

$$|u_i^n(x) - u_i^0(x)|^2 \leq (1+\varepsilon)|u_i^n(x) - \tilde{u}_i^n(x)|^2 + C_\varepsilon|\tilde{u}_i^n(x) - u(t_i, x)|^2$$
$$\leq (1+\varepsilon)[G(u^n - \tilde{u}^n)|x|^2 + H(u^n - \tilde{u}^n)] + C_\varepsilon(1+|x|^2)h$$
$$\leq [(1+\varepsilon)c_2(\lambda_1, h, \bar{L}, \bar{G})G(u^n - u^0) + C_\varepsilon h]|x|^2$$
$$\quad + (1+\varepsilon)c_2(\lambda_1, h, \bar{L}, \bar{G})$$
$$\quad \times [H(u^n - u^0) + [B_1 + A_2\bar{H}]TG(u^n - u^0)] + C_\varepsilon h.$$

Now for any $c_2(L_1, L_1) < c_2 < 1$, we can choose $\bar{L}, \bar{G}$ and $\varepsilon$ appropriately such that, for $h$ small enough,

$$(1+\varepsilon)c_2(\lambda_1, h, \bar{L}, \bar{G}) \leq c_2.$$

Then we get

$$|u_i^n(x) - u_i^0(x)|^2 \leq [c_2 G(u^n - u^0) + C_\varepsilon h]|x|^2 \tag{6.5}$$
$$+ c_2 H(u^n - u^0) + C_\varepsilon G(u^n - u^0) + C_\varepsilon h.$$

We now follow the arguments in the proof of Theorem 5.1(iii). Fix some $G_0, H_0$ such that

$$|u_i^n(x) - u_i^0(x)|^2 \leq G_0|x|^2 + H_0.$$

For $\nu = 1, 2, \ldots$, denote

$$G_\nu \triangleq c_2 G_{\nu-1} + C_\varepsilon h, \qquad H_\nu \triangleq c_2 H_{\nu-1} + C_\varepsilon G_{\nu-1} + C_\varepsilon h.$$

Then (6.5) implies that

$$|u_i^n(x) - u_i^0(x)|^2 \leq G_\nu |x|^2 + H_\nu \qquad \forall \nu. \tag{6.6}$$

Note that

$$G_\nu = G_0 c_2^\nu + C_\varepsilon h \frac{1 - c_2^\nu}{1 - c_2},$$

$$H_\nu = H_0 c_2^\nu + C_\varepsilon G_0 \nu c_2^\nu + \frac{C_\varepsilon h}{1 - c_2}\left[\frac{1 - c_2^\nu}{1 - c_2} - \nu c_2^\nu\right] + C_\varepsilon h \frac{1 - c_2^\nu}{1 - c_2}.$$

Since $c_2 < 1$, and sending $\nu \to \infty$ in (6.6), we get

$$|u_i^n(x) - u_i^0(x)|^2 \leq \frac{C_\varepsilon h}{1 - c_2}|x|^2 + \frac{C_\varepsilon h}{(1 - c_2)^2}.$$

The proof is complete. □

As a direct consequence of Theorems 5.1 and 6.3, we have the following:

THEOREM 6.4. *Under the assumptions of Theorem 6.3 we have, for any $c_2(L_1, L_1) < c_2 < 1$ and for h small enough,*

$$|u_i^{n,m}(x) - u(t_i, x)|^2 \leq C(1 + |x|^2)[mc_2^m + h].$$

We close the theoretical part of this paper with a precise version of the generic Theorem 2.3.

THEOREM 6.5. *Under the assumptions of Theorem 6.3 we have, for any $c_2(L_1, L_1) < c_2 < 1$ and for h small enough,*

$$\sup_{1 \leq i \leq n} E\left\{\sup_{t \in [t_{i-1}, t_i]}[|X_t - X_{i-1}^{n,m}|^2 + |Y_t - Y_{i-1}^{n,m}|^2]\right\}$$
$$+ \sum_{i=1}^n E\left\{\int_{t_{i-1}}^{t_i} |Z_t - \hat{Z}_{i-1}^{n,m}|^2 \, dt\right\} \leq C(1 + |x|^2)[mc_2^m + h].$$

PROOF. By (6.4), it suffices to prove

$$(6.7) \quad \sup_{0 \leq i \leq n} E\{|\Delta X_i|^2 + |\Delta Y_i|^2\} + h \sum_{i=0}^{n-1} E\{|\Delta Z_i|^2\} \leq C(1+|x|^2)[mc_2^m + h],$$

where

$$\Delta X_i \triangleq \tilde{X}_i^n - X_i^{n,m}, \qquad \Delta Y_i \triangleq \tilde{Y}_i^n - Y_i^{n,m}, \qquad \Delta Z_i \triangleq \tilde{Z}_i^n - \hat{Z}_i^{n,m}.$$

First, by Lemma 3.2 with $\lambda_1 = 1$, we get

$$E\{|\Delta X_{i+1}|^2\} \leq E\{(1+Ch)|\Delta X_i|^2 + Ch|u(t_i, X_i^{n,m}) - u_i^{n,m}(X_i^{n,m})|^2\}$$
$$\leq (1+Ch)E\{|\Delta X_i|^2\} + C(1+|x|^2)[mc_2^m + h]h.$$

Since $\Delta X_0 = 0$, we have

$$\sum_{0 \leq i \leq n} E\{|\Delta X_i|^2\} \leq C(1+|x|^2)[mc_2^m + h].$$

Next, choose $\lambda_2 = \lambda_3 = \frac{1}{5}$ and $h$ small enough so that $A_3 \leq \frac{1}{2}$. Applying Lemma 3.3, we obtain

$$E\{|\Delta Y_i|^2 + \tfrac{1}{2}h|\Delta Z_i|^2\} \leq E\{(1+Ch)|\Delta Y_{i+1}|^2 + Ch|\Delta X_i|^2\}.$$

Since

$$|\Delta Y_n|^2 = |g(\tilde{X}_n^n) - g(X_n^{n,m})|^2 \leq C|\Delta X_n|^2,$$

we can easily get

$$\sup_{0 \leq i \leq n} E\{|\Delta Y_i|^2\} + h \sum_{i=0}^{n-1} E\{|\Delta Z_i|^2\}$$
$$\leq C \sup_{0 \leq i \leq n} E\{|\Delta X_i|^2\} \leq C(1+|x|^2)[mc_2^m + h].$$

This proves (6.7) and hence the theorem. □

**7. A numerical algorithm.** We now briefly explain how the discretized Markovian iteration above can be transformed into a numerical algorithm which is viable also for high-dimensional problems. To this end, we replace the conditional expectations by a simulation based least squares regression estimator, as was suggested, for example, by Gobet, Lemor and Warin [14] and Bender and Denk [3] in the context of decoupled FBSDEs. An alternative estimator based on Malliavin calculus is discussed in Bouchard and Touzi [4] for decoupled FBSDEs. A quantization algorithm for reflected BSDEs is presented in Bally and Pagès [2].

For the reader's convenience, we spell out our algorithm for the coupled case. While a convergence analysis is out of the scope of the present paper, we will illustrate the algorithm by some numerical examples in the next section.

We assume that the number of time steps $n$ is fixed for the remainder of this section. In the algorithm conditional expectations are first replaced by orthogonal projections on $K$ basis functions. Then the orthogonal projections are approximated by simulating $\Lambda$ trajectories. Hence, the algorithm can be described for the one-dimensional case iteratively as follows. It is straightforward how this extends to the multi-dimensional case:

- Fix some $x_0$. Set $\bar{u}_i^{n,0,K,\Lambda}(x) \triangleq 0$.
- Sample $\Lambda$ independent copies of the time discretized Brownian motion $W_{t_i}^\lambda$, $i = 0, \ldots, n$, $\lambda = 1, \ldots, \Lambda$, starting in 0 and denote the corresponding increments by $\Delta W_i^\lambda$.
- Suppose $\bar{u}_i^{n,m-1,K,\Lambda}(x)$ is already constructed. Let $\bar{X}_0^{n,m,\lambda} \triangleq x_0$ and

$$\bar{X}_{i+1}^{n,m,\lambda} \triangleq \bar{X}_i^{n,m,\lambda} + b(t_i, \bar{X}_i^{n,m,\lambda}, \bar{u}_i^{n,m-1,K,\Lambda}(\bar{X}_i^{n,m,\lambda}))h$$
$$+ \sigma(t_i, \bar{X}_i^{n,m,\lambda}, \bar{u}_i^{n,m-1,K,\Lambda}(\bar{X}_i^{n,m,\lambda}))\Delta W_{i+1}^\lambda,$$

where—for notational convenience—we suppress the dependence of $\bar{X}_i^{n,m,\lambda}$ on $K$ through $\bar{u}^{n,m-1,K,\Lambda}$. Note, $\bar{X}_i^{n,m,\lambda_0}$ depends on all Brownian increments $\Delta W_i^\lambda$, $i = 1, \ldots, n$, $\lambda = 1, \ldots, \Lambda$, through $\bar{u}_i^{n,m-1,K,\Lambda}$. While we expect that this dependence will make a convergence analysis difficult, the examples below indicate that the algorithm works without re-simulating the Brownian paths in every iteration step.

- Choose a set of Lipschitz continuous basis functions

$$\mathcal{B}_i^{n,m,K} \triangleq \{\eta_i^{n,m,k}(x),\ k = 1, \ldots, K\}$$

such that

(7.1) $$\{\eta_i^{n,m,k}(\bar{X}_i^{n,m,\lambda}),\ k = 1, \ldots, K\}$$

forms a subset of $L^2(\Omega)$. From the construction below, it will become evident that $\bar{u}_i^{n,m,K,\Lambda}(x)$ inherits the Lipschitz continuity from the basis functions. This feature seems to be important to ensure that the discretized forward equations for $\bar{X}^{n,m+1,\lambda}$ do not explode.

- Define, for $i = n - 1, \ldots, 1$,

$$\bar{u}_n^{n,m,K,\Lambda}(x) \triangleq g(x), \qquad \bar{v}_n^{n,m,K,\Lambda}(x) \triangleq 0,$$
$$\bar{Y}_{i+1}^{n,m,K,\Lambda} \triangleq \bar{u}_{i+1}^{n,m,K,\Lambda}(\bar{X}_{i+1}^{n,m,\lambda}),$$
$$\bar{v}_i^{n,m,K,\Lambda}(x) \triangleq \arg\inf\left\{\frac{1}{\Lambda}\sum_{\lambda=1}^{\Lambda}\left|\frac{1}{h}\bar{Y}_{i+1}^{n,m,K,\Lambda}\Delta W_{i+1}^\lambda - V(\bar{X}_i^{n,m,\lambda})\right|^2\right.;$$

$$V \in \mathrm{span}(\mathcal{B}_i^{n,m,K})\bigg\},$$

$$\bar{Z}_i^{n,m,K,\lambda} \triangleq \bar{v}_i^{n,m,K,\Lambda}(\bar{X}_i^{n,m,\lambda}),$$

$$\bar{u}_i^{n,m,K,\Lambda}(x) \triangleq \arg\inf\bigg\{\frac{1}{\Lambda}\sum_{\lambda=1}^{\Lambda}|f(t_i,\bar{X}_i^{n,m,\lambda},\bar{Y}_{i+1}^{n,m,K,\lambda},\bar{Z}_i^{n,m,K,\lambda})h$$

$$+ \bar{Y}_{i+1}^{n,m,K,\lambda} - U(\bar{X}_i^{n,m,\lambda})|^2;$$

$$U \in \mathrm{span}(\mathcal{B}_i^{n,m,K})\bigg\}.$$

Note that the minimization problems are linear least squares problems, which can be easily implemented.

- Let

$$\bar{Y}_1^{n,m,K,\lambda} \triangleq \bar{u}_1^{n,m,K,\Lambda}(\bar{X}_1^{n,m,\lambda}),$$

$$\bar{Z}_0^{n,m,K,\lambda} \triangleq \frac{1}{\Lambda}\sum_{\bar{\lambda}=1}^{\Lambda}\frac{1}{h}\bar{Y}_1^{n,m,K,\bar{\lambda}}\Delta W_1^{\bar{\lambda}},$$

$$\bar{Y}_0^{n,m,K,\lambda} \triangleq \frac{1}{\Lambda}\sum_{\bar{\lambda}=1}^{\Lambda}\bar{Y}_1^{n,m,K,\bar{\lambda}} + f(0,x_0,\bar{Y}_1^{n,m,K,\bar{\lambda}},\bar{Z}_0^{n,m,K,\bar{\lambda}})h.$$

We expect that the thus constructed $(\bar{X}^{n,m,\lambda},\bar{Y}^{n,m,K,\lambda},\bar{Z}^{n,m,K,\lambda})$ are "close" to $(X^{n,m,\lambda},Y^{n,m,\lambda},\hat{Z}^{n,m,\lambda})$, the solution of the discretized Markovian iteration (1.8) with the Brownian motion $W$ replaced by $W^\lambda$, if the basis functions are chosen appropriately and the number $\Lambda$ of simulated paths is sufficiently large. While an analysis of the error by estimating the conditional expectations is left to future research, the numerical examples in the next section support this conjecture.

**8. Numerical examples.** For the simulations, we consider the example

$$\begin{cases} X_{d,t} = x_{d,0} + \int_0^t \sigma Y_u\, dW_{d,u}, \\ Y_t = \sum_{d=1}^D \sin(X_{d,T}) + \int_t^T -rY_u + \tfrac{1}{2}e^{-3r(T-u)}\sigma^2\left(\sum_{d=1}^D \sin(X_{d,u})\right)^3 du \\ \qquad - \int_t^T \sum_{d=1}^D Z_{d,u}\, dW_{d,t}, \end{cases}$$

where $W_{d,t}, d=1,\ldots,D$, is a $D$-dimensional Brownian motion and $\sigma > 0$, $r$, $x_{d,0}$ are constants. Note that the corresponding differential operator degenerates at $y=0$.

By Itô's formula, one can easily check that this FBSDE decouples via the relation

$$Y_t = e^{-r(T-t)} \sum_{d=1}^{D} \sin(X_{d,t}). \tag{8.1}$$

Note that for small $\sigma$ the weak coupling condition of $Y$ into $X$ is satisfied, while, for large $\sigma$, the monotonicity condition of $f$ can be fulfilled by choosing $r$ large enough.

In the simulations we replace conditional expectations by least squares regression as explained above with the "canonical" basis functions

$$1, \qquad x_d, \qquad 1 \leq d \leq D, \qquad (-R) \vee (x_d x_q) \wedge R, \qquad 1 \leq d \leq q \leq D,$$

that is, monomials up to order two in $x = (x_1, \ldots, x_D)$. The truncation constant $R$ guarantees that the basis functions are Lipschitz continuous. We set

$$R \triangleq 10, \qquad X_{d,0} \triangleq \frac{\pi}{2}, \qquad 1 \leq d \leq D,$$

$$T \triangleq 1, \qquad \Lambda \triangleq 50000, \qquad n \triangleq 50,$$

unless otherwise stated. With this initial condition, we get $Y_0 = De^{-r(T-t)}$. Recall also that the estimator $\bar{Y}_0^{n,m,K,\lambda}$ of $Y_0$ does not depend on $\lambda$ and is denoted by $\bar{Y}_0^{n,m,K}$ from now on.

Figure 1 illustrates the convergence of the iteration in the case of a four-dimensional state space ($D = 4$). Both figures display the absolute error $|\bar{Y}_0^{n,m,K} - Y_0|$ as a function of the number of iterations $m$. In Figure 1(a) the case $r = 0$ (no monotonicity) is considered for several values of $\sigma$ which represent different influences of the coupling. In Figure 1(b) the coupling parameter $\sigma = 0.4$ is fixed, while the strength of the monotonicity varies by different values of $r$. In general, we observe that the iteration converges extremely fast, as could be expected in view of Theorem 2.1 which states $mc^m$, for some $c$, $0 < c < 1$, as rate of convergence. From the proof of this theorem we know that $c$ is the smaller, the weaker the coupling or the stronger the monotonicity is. This explains the faster convergence observed for small values of $\sigma$ and large values of $r$.

The influence of the time partition is displayed in Figure 2. It shows the absolute error $|\bar{Y}_0^{n,m_{\text{stop}},K} - Y_0|$ as a function of the number of time points $n$. We stop the iteration when two consecutive estimates $\bar{Y}_0^{n,m,K}$ are within a distance of $10^{-4}$. This iteration level is denoted $m_{\text{stop}}$. The observed convergence rate is in accordance with $1/\sqrt{n}$ as derived in Theorem 2.2.

Finally, we demonstrate that the space dimension four is no limitation for the proposed algorithm. To this end, we consider the 10-dimensional case with the parameter values $r = 0$, $\sigma = 0.1$ in Figure 3. Under the same stopping criterion as above, the iteration terminates after 12 steps. Recall that, from (8.1), we can also approximate the true value of $X$ via the usual Euler scheme (applying the same simulated Brownian increments $\Delta W_i^\lambda$). The corresponding approximation along the $\lambda$th path is denoted $\check{X}_i^{n,\lambda}$, and hence,

$$\check{Y}_i^{n,\lambda} = e^{-r(T-t_i)} \sum_{d=1}^{D} \sin(\check{X}_{d,i}^{n,\lambda})$$

may be considered a close approximation of $Y_{t_i}$. In Figure 3 we display, for the 10-dimensional case, a comparison between a typical path of $\check{Y}_i^{n,\lambda_0}$ (dashed line) and $\bar{Y}_i^{n,m_{\text{stop}},K,\lambda_0}$ (solid line), as well as the (absolute) empirical mean square error between $\check{Y}_i^{n,\lambda}$ and $\bar{Y}_i^{n,m_{\text{stop}},K,\lambda}$, $\lambda = 1,\ldots,\Lambda$. Precisely, Figure 3(b) shows

$$\frac{1}{\Lambda} \sum_{\lambda=1}^{\Lambda} |\bar{Y}_i^{n,m_{\text{stop}},K,\lambda} - \check{Y}_i^{n,\lambda}|^2$$

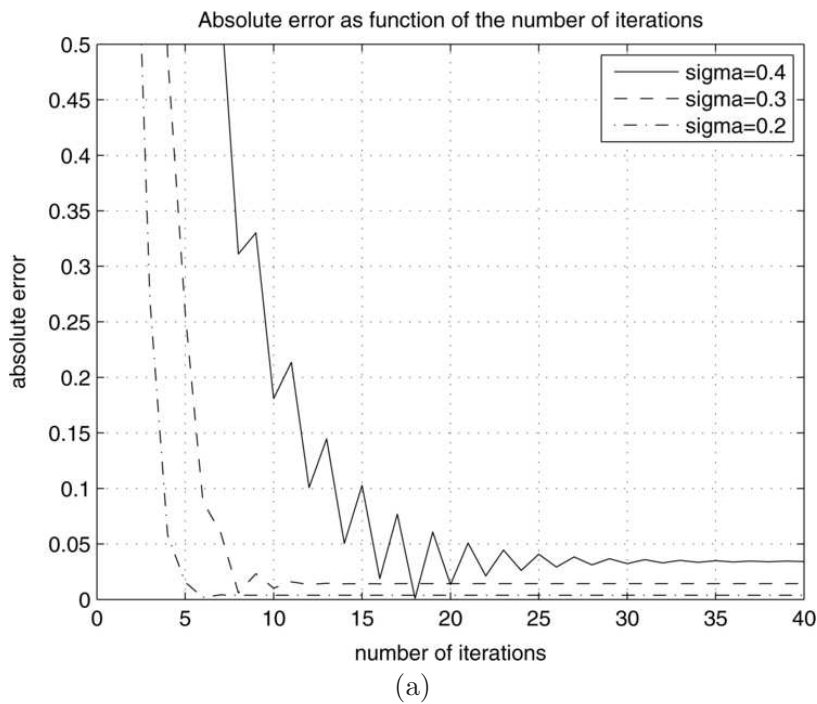

(a)

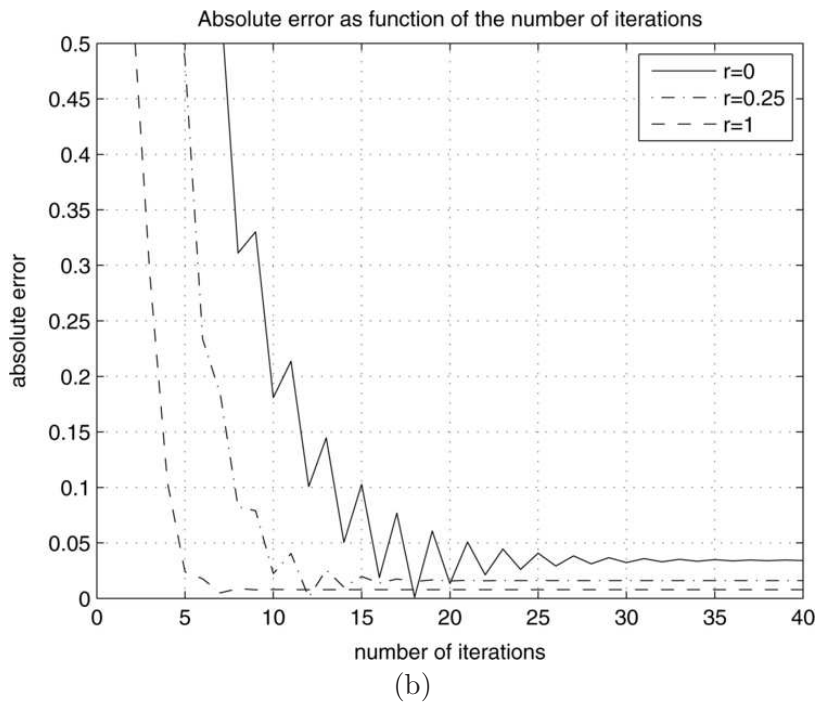

(b)

FIG. 1. *Convergence of the iteration for different choices of $\sigma$ and $r$.*

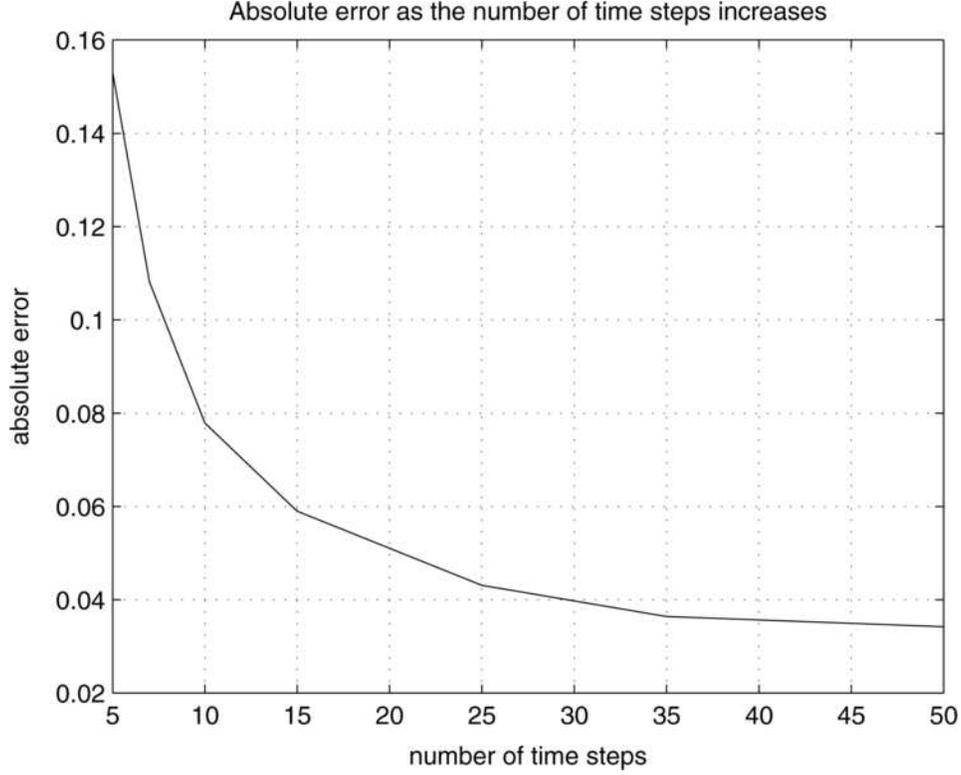

Fig. 2. *Convergence as the time partition becomes finer for $r = 0$, $\sigma = 0.4$, $D = 4$.*

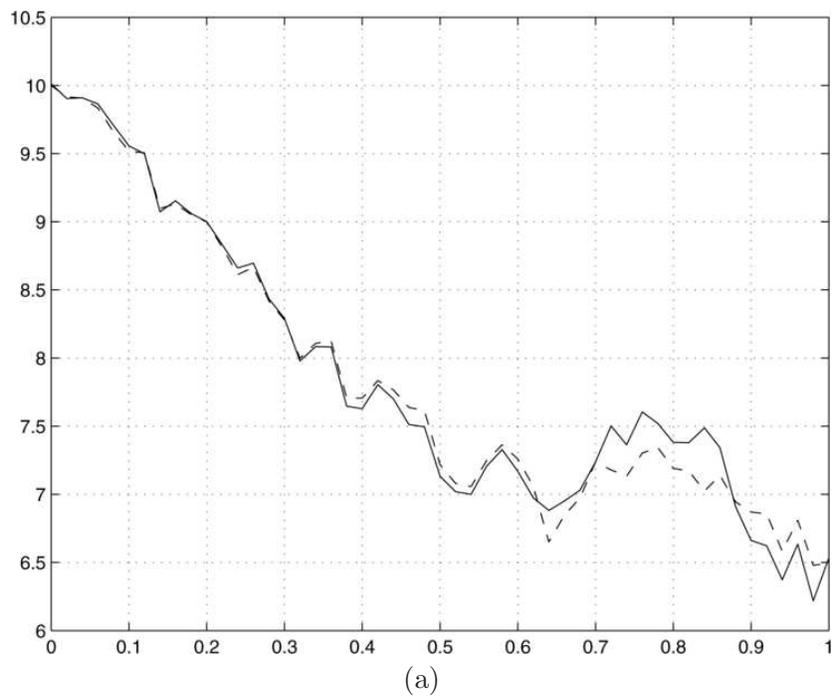

(a)

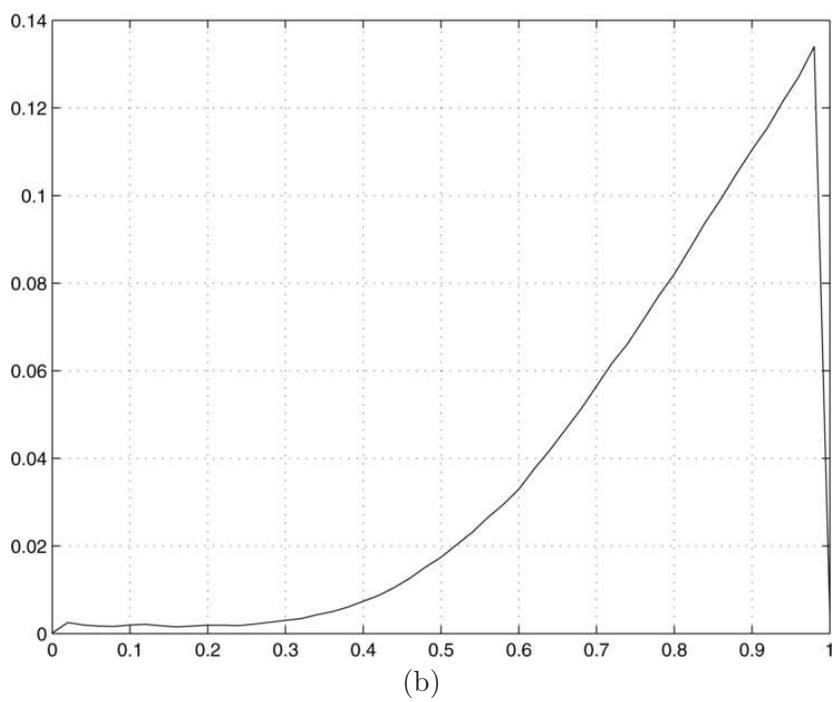

(b)

Fig. 3. *Typical path of $\bar{Y}_i^{n,m_{\mathrm{stop}},K,\lambda}$ and mean square error for $\sigma = 0.1$, $r = 0$, $D = 10$.*

as function of time.

REMARK 8.1. Figure 3 shows a larger mean square error close to terminal time than close to initial time. This is a rather typical feature, when conditional expectations are estimated by the above least squares method. It can be explained by interpreting this method as creating a stochastic mesh (see Glasserman [13]) and observing that this mesh is typically much finer close to initial time than close to terminal time. Hence, the error close to terminal time, observed in Figure 3, cannot be diminished by solely increasing the number of iterations, but by improving the quality of the conditional expectation estimator. A generic trick to improve the quality close to terminal time is to add the terminal function $g$ to the basis.

**Acknowledgment.** 

INSTITUTE FOR MATHEMATICAL STOCHASTICS
TU BRAUNSCHWEIG
POCKELSSTR. 14
D-38106 BRAUNSCHWEIG
GERMANY
E-MAIL: c.bender@tu-bs.de

DEPARTMENT OF MATHEMATICS
UNIVERSITY OF SOUTHERN CALIFORNIA
LOS ANGELES, CALIFORNIA 90089
USA
E-MAIL: jianfenz@usc.edu